\documentstyle[12pt]{article}
\def\mineappendix{
        \setcounter{section}{1}
        \setcounter{subsection}{0}
        \def\thesection{\Alph{section}}
        \def\sectionap{\@startsection  {section}{1}{\z@}
                        {-3.5ex plus-1ex minus-.2ex} {0ex plus.2ex}
                        {\reset@font\Large\bf  Appendix:  \, }
                        }
        }
\makeatother
\def\Proclaim #1. #2\par{\bigbreak\noindent{\sc#1.\enspace}{\it#2}\par}

\font\Bbbfont=msbm10
\newfam\msbfam
\textfont\msbfam=\Bbbfont \textfont\msbfam=\Bbbfont

\def\pa{\partial}
\def\al{\alpha}

\def\si{\sigma}

\def\wt{\widetilde}

\def\wt{\widetilde}
\def\ch{\mbox{ch}}
\def\sh{\mbox{sh}}

\def\nn{\nonumber}
\def\Ga{\Gamma}

\def\om{\omega}

\def\nd{\noindent}

\newtheorem{Theorem}{Theorem}%[section]
\newtheorem{Lemma}{Lemma}%[section]

\newtheorem{Remark}{Remark}

\title{ On the $1+2$ Dimensional Isotropic Landau-Lifshitz Equation}
\author{Qing Ding\\
School of Mathematical Sciences, Fudan University\\ Shanghai
200433, China\\
 E-mail address: qding@fudan.edu.cn}
\date{}

\begin{document}
\maketitle

\begin{abstract}
By using the geometric concept of PDEs with prescribed curvature
representation, we prove that the $1+2$ dimensional isotropic
Landau-Lifshitz equation
%, or in other words, the Schr\"odinger
%map from ${\bf R}^2$ to $S^2\hookrightarrow {\bf R}^3$,
is gauge
equivalent to a $1+2$ dimensional nonlinear Schr\"odinger-type
system. From the nonlinear Schr\"odinger-type system, we construct
blowing up $H^3({\bf R}^2)$-solutions to the Landau-Lifshitz
equation, which reveals the blow up phenomenon of the equation.
\end{abstract}

\bigskip

\section *{\S 1. Introduction}
The isotropic Landau-Lifshitz equations, or in other words, the
generalized Heisenberg models for a continuous ferromagnetic spin
%system with classical spin
vector $s=(s_1,s_2,s_3)$ $\in
S^2\hookrightarrow {\bf R}^3$ (see, for example,
\cite{LaL,PaT,CSU,SuSu}),
\begin{eqnarray}
{\bf s}_t={\bf s}{\times}\Delta_{{\bf R}^n}{\bf  s}, \quad x\in
{\bf R}^n,\quad n=1,2,3,\cdots \label{0}
\end{eqnarray}
are important equations in spin magnetic fields in physics. These
equations exhibit a rich variety of dynamical properties of a spin
vector in different backgrounds.

Though Eq.(\ref{0}) have a unified version of expressions for
$n\ge1$, there are great differences between dynamical properties
of Eq.(\ref{0})  with $n=1$ and those of Eq.(\ref{0}) with
$n\ge2$. When $n=1$, Eq.(\ref{0}) is integrable and it can be
solved by the method of inverse scattering techniques
(\cite{FaT,bib:Gu}). Furthermore, Eq.(\ref{0}) with $n=1$ is gauge
equivalent to the nonlinear Schr\"odinger equation of attractive
type: $i\phi_t+\phi_{xx}+2 |\phi|^2\phi=0$ (\cite{ZaT}) (its dual
version was proved in \cite{d1,d3}).
% and, meanwhile, Eq.(\ref{1}) is gauge equivalent to the nonlinear
%Schr\"odinger equation of repulsive type: $i\phi_t+\phi_{xx}-2
%|\phi|^2\phi=0$ (\cite{bib:dq1}).
%From this one can obtains the
%global existence of the Cauchy problem of Eq.(\ref{0}) (resp.
%Eq.(\ref{1})) directly from those of nonlinear Schr\"odinger
%equations.
When $n\ge2$, Eq.(\ref{0}) are non-integrable and the
understanding of their dynamical properties becomes much
difficulty. In 1975,  Belavin and Polyakov (\cite{bib:BP}) were
the first who paid attention to the construction of topological
static solutions (or in other words, Belavin-Polyakov instantons)
to Eq.(\ref{0}) with $n=2$, i.e., the 1+2 dimensional isotropic
Landau-Lishitz equation, by applying the technique of the
stereographic projection $S^2\to C$. These solutions are used for
the description of domain walls, magnetic bubbles and the
metastatic states in ferromagnetic spin fields. The consequences
of such a study are referred to \cite{PaT} and references therein.
In 1986, Sulem, Sulem and Bardos proved in \cite{SuSu} the local
$H^{m+1}({\bf R }^n)$-existence of solutions to the Cauchy problem
of Eq.(\ref{0}) by the difference method and also the global
$W^{m+1,6}({\bf R}^n)$-existence if initial data is small enough.
Very recently, in \cite{NSU} Nahmod, Stefanov and Uhlenbeck proved
the local well-posedness of the Cauchy problem of the 1+2
dimensional Landau-Lifshitz equation and its dual equation by
applying their equivalent equations so-called the modified
Schr\"odinger map equations.

The isotropic Landau-Lifshitz equations (\ref{0}) are special
cases of so-called Schr\"odinger maps (\cite{CSU,GS,NSU}) or
Schr\"odinger flows (\cite{FMP,DW2,TK,d1}) in geometry. The
Schr\"odinger  map from a Riemannian manifold $(M,g)$ to a
K\"ahler manifold $(N,\omega)$ is defined to be the (infinite
dimensional) Hamiltonian system of the energy function
$E(u)=\int_M|\nabla u|^2dv_g$ on the mapping space $C^k(M,N)$ for
some $k>0$. More explicitly, let $J$ be a compatible complex
structure on $N$ such that $h(\cdot,\cdot)=\om(\cdot,J\cdot)$ is a
Riemannian metric on $N$ and  we denote by $\nabla E(u)$ the
gradient of the function $E(u)$ with respect to the inner product
$ <v,w>_u=\int_{M}h(u)(v,w),\quad \forall u\in X, \forall v,w\in
T_uX \label{1}$, on $X$, then the corresponding Hamiltonian vector
field $V_{E(u)}$ can be expressed explicitly as
$V_{E(u)}=J(u)\nabla E(u)$. Thus the Schr\"odinger map from
$(M,g)$ into $(N,\om)$ is represented by
\begin{eqnarray}
u_t=J(u)\nabla E(u).\nn
\end{eqnarray}
It is easy to verify that the gradient $\nabla E(u)$ is exactly
the tension field $\tau(u)$ of map $u:M\to N$. In a local
coordinates
$$
\tau^l(u)=\Delta_{M}u^l+\Ga^l_{jk}(u){{\pa u^j} \over{\pa
x^{\al}}}{{\pa u^k}\over{\pa x^{\beta}}}g^{\al\beta}, \quad 1\le
l\le \hbox{dim}N,
$$
where $\Delta_M$ denotes the Laplacian operator on $(M,g)$ with
the given metric $g=(g_{\alpha\beta})$, $(g^{\alpha\beta})$ is the
inverse matrix of $(g_{\alpha\beta})$ and $\Gamma^i_{jk}$ are the
Christoffel symbols of the target manifold $(N,h)$. So the
Schr\"odinger  map from $M$ into $N$ can also be written as:
\begin{eqnarray}
u_t=J(u)\tau(u). \nn
\end{eqnarray}
It is a straightforward verification that the Schr\"odinger map
from an Euclidean $n$-space ${\bf R}^n$ to the 2-sphere
$S^2\hookrightarrow{\bf R}^3$ is exactly the isotropic
Landau-Lifshitz equation (\ref{0}) (for example see \cite{CSU}).

The main object in the study of Schr\"odinger maps is, of course,
the solvability of the corresponding Cauchy or initial-boundary
value problem and its solutions' behaviors. However, comparing to
those of heat flows or wave maps \cite{Str}, this study is still
at the beginning stage. Except the results stated above, the
following results dealing with Schr\"odinger maps should also be
mentioned. In 1999, Terng and Uhlenbeck showed the global
existence of the Cauchy problem of Schr\"odinger maps from ${\bf
R}^1$ to complex compact Grassmannians in \cite{TK}. Chang, Shatah
and Uhlenbeck proved in \cite{CSU} the global existence and
uniqueness of smooth solution to the Cauchy problem of
Schr\"odinger maps from $R^1$ to compact Riemainnian surfaces and
also the $W^{2,4}({\bf R}^2)$-global existence result of the
radial Schr\"odinger maps from ${\bf R}^2$ to an arbitrary compact
Riemann surface for small initial data. W.Y. Ding and Wang proved
in \cite{DW2} the existence of local smooth or global weak
solutions to the Cauchy problem of Schr\"odinger maps from a
compact Riemannian manifold or Euclidean space ${\bf R}^n$ to a
compact K\"ahler manifold. Grillakis and Stefanopoulos displayed
conservation laws and localized energy estimates of Schr\"odinger
maps to Riemannian surfaces. However, it is widely believed that a
Schr\"odinger map with dimensions of the starting manifold are
greater than 1 may develop singularities in finite time in
general. This fundamental problem in the study of Schr\"odinger
maps was proposed by W.Y. Ding as a unsolved question in
\cite{Di}. In fact, the same question for the higher dimensional
Landau-Lifshitz equations has been existed for a long time. Blow
up corresponds to the self-trapping and intense focussing
phenomena of classical ferromagnetic spin fields.

In this paper, we first display that the $1+2$ dimensional
iostropic Landau-Lifshitz equation: ${\bf s}_t={\bf
s}{\times}\left({{\partial^2}\over{\partial x^2}}+
{{\partial^2}\over{\partial y^2}}\right){\bf s}$ is gauge
equivalent to the following $1+2$ dimensional nonlinear
Schr\"odinger-type system
\begin{eqnarray}
\left\{\begin{array}{c} iq_t-q_{z{\bar z}}+2u
q-2({\bar p}q)_z+2pq_{\bar z}+4|p|^2q=0\\
ir_t+r_{z{\bar z}}-2u r-2({\bar p}r)_z+2pr_{\bar z}-4|p|^2r=0\\
ip_t=(qr)_{\bar z}-u_z\qquad\qquad\qquad\qquad\qquad\qquad
\end{array} \right.    \nn
\end{eqnarray}
with the additional restrictions: ${\bar p}_z+p_{\bar z}=
|q|^2-|r|^2,\quad {\bar r}_z+q_{\bar z}=2(p{\bar r}-{\bar p}q)$ by
using the geometric concept of gauge equivalence for PDEs with
prescribed curvature representation developed in \cite {d4,dz}.
This nonlinear Schr\"odinger-type system is very different from
the modified Schr\"odinger map equation obtained by Nahmod,
Stefanov and Uhlenbeck in \cite{NSU}. Then, by characterizing some
analytic properties of the above nonlinear Schr\"odinger-type
system, we show the existence of blowing up $H^{3}({\bf
R}^2)$-solutions to the $1+2$ dimensional Landau-Lifshitz equation
if the initial data is chosen to be technically small. This
blowing up result does not contradict to the global existence of
$W^{m+1,6}({\bf R}^2)$-solutions due to Sulem, Sulem and Bardos in
\cite{SuSu}. On the contrary, it reflects some new interesting and
mysterious properties of the Landau-Lifshitz equation.

This paper is organized as follows. In the section 2, we  shall
transfer the $1+2$ dimensional isotropic Landau-Lifshitz equation
to its gauge equivalent nonlinear Schr\"odinger-type system. In
section 3, by applying the nonlinear Schr\"odinger-type system, we
construct blowing up $H^{3}({\bf R}^2)$-solutions to the $1+2$
dimensional isotropic Landau-Lifshitz equation  and give some
remarks.

\section * {\S 2. Gauge Equivalence}
In almost known results dealing with the  Landau-Lifshitz
equations and Schr\"odinger maps, the first step of the study is
to transform the original equation to a nonlinear
Schr\"odinger-type equation. In this section we shall apply the
geometric concept of gauge transformations to transfer the 1+2
dimensional Landau-Lifshitz equation into a nonlinear
Schr\"odinger-type system in the category of (nonzero) prescribed
curvature formulations. The zero curvature formulation in
integrable theory is a main indication of integrability of a
soliton equation. Zarkharov and Takhtajan introduced in \cite{ZaT}
the geometric concept of gauge equivalence between two soliton
equations which provides a useful tool in the study of integrable
equations. In \cite{d4,dz} the author and his collaborator found
that the geometric concept of gauge equivalence can be generalized
to differential equations with prescribed curvature representation
and then displayed the gauge equivalent structures of the  $1+1$
dimensional anisotropic Landau-Lifshitz equation (which give an
affirmative answer to a question proposed in \cite{BoY}) and the
modified nonlinear Schr\"odinger equation. Now we find that it is
also applicable to the present $1+2$ dimensional isotropic
Landau-Lifshitz equation.

Let us explicitly write down the $1+2$ dimensional isotropic
Landau-Lifshitz equation  Eq.({\ref{0}) as follows:
\begin{eqnarray}
{\bf s}_t={\bf s}{\times}\left({{\partial^2}\over{\partial x^2}}+
{{\partial^2}\over{\partial y^2}}\right){\bf s}, \label{4}
\end{eqnarray}
where $(x,y)$ is the standard coordinates of the Euclidean plane
${\bf R}^2$. We convert it into the matrix form,
\begin{eqnarray}
S_t=-{1\over2i}[S, S_{z{\bar z}}], \label{5}
\end{eqnarray}
where  $S=\left(\begin{array}{cc}
       s_3 &s _1+is_2\\
       s_1-is_2 & -s_3
       \end{array}\right)$
        with $S^2=I$ ($I$ denotes the unit matrix as usual) and
         $z={{x+iy}\over2}$, ${\bar z}={{x-iy}\over2}$ which are
        ${1\over2}$ scaling of the usual complex versions of the real
        variables $x$ and $y$.
In order to present Eq.(\ref{5}) as an equation with prescribed
curvature representation,  let's set
\begin{eqnarray}
A= Vd{\bar z} -i\lambda S d{z} +\lambda(2V +SS_{\bar z}+2\al S)dt
\label{6}
\end{eqnarray}
and
\begin{eqnarray}
K&=&\Bigg\{-V_t+\lambda (2V +SS_{\bar z}+2\al S)_{\bar
z}-[V,\lambda(2V
+SS_{\bar z}+2\al S)] \Bigg\}d{\bar z} \wedge dt \nonumber \\
&&+\lambda\Bigg({1\over2}[S_{z},S_{\bar z}]+2(\al S)_{z}\Bigg)
d{z} \wedge dt, \label{K}
\end{eqnarray}
where $\lambda$ is a spectral parameter which is independent of
$t$,$z$ and $\bar z$, $V=V(\lambda,z,\bar z,t)$ is a
$2\times2$-matrix satisfying the equation:
\begin{eqnarray}
 (-i\lambda S)_{\bar z}-V_{z} +[-i\lambda S,V]=0 \label{V}
\end{eqnarray}
and $\al={1\over{2}}\hbox{tr}(G^{-1}G_{\bar z}\sigma_3)$ is a
scalar function, where $G$ is an $SU(2)$-matrix solving (\ref{14})
below.
%(note that, by a straightforward computation,
%$\hbox{tr}(SS_{z}S_{\bar z})=  2i(s_1({s_2}_{\bar
%z}{s_3}_z-{s_2}_z{s_3}_{\bar z})+s_2({s_3}_{\bar
%z}{s_1}_z-{s_3}_z{s_1}_{\bar z})+s_3({s_1}_{\bar
%z}{s_2}_z-{s_1}_z{s_2}_{\bar z}))$ is real in this case).
$d+A$ can be geometrically interpreted as defining a connection on
a trivial $SU(2)$-principal bundle  over ${\bf R}^3$ (the space of
the independent variables $x$, $y$ and $t$). Then it is a
straightforward computation that the curvature $F_A$ of the
connection $d+A$ is
\begin{eqnarray}
F_A&=&dA-A\wedge A\nn\\
 &=&\Bigg\{-V_t+\lambda (2V +SS_{\bar
z}+2\al S)_{\bar z}-[V,\lambda(2V
+SS_{\bar z}+2\al S)] \Bigg\}d{\bar z} \wedge dt \nonumber \\
&&+\Bigg\{\lambda\bigg(iS_t+(SS_{\bar z}+2\al S)_{z}\bigg)\Bigg\}
d{z} \wedge dt \nonumber
\end{eqnarray}
and hence Eq.(\ref{5}) is equivalent to holding the following
prescribed curvature condition:
\begin{eqnarray}
F_A=dA-A\wedge A=K. \label{7}
\end{eqnarray}
We would like to point out that, though we cannot give an explicit
solution $V$ to Eq.(\ref{V}), it still OK for us to prove our
desired conclusion, as we shall see below.

It is well-known that, in the Yang-Mills gauge theory, there are
gauge transformations $A\to {\hat A}=dG G^{-1}+ G A G^{-1}$, for
 $G\in C^{\infty}({\bf R}^3, SU(2))$
such that $F_{\hat A}=G F_A G^{-1}$ under the gauge
transformation.

\begin{Theorem} There is a gauge function $G(t,z,{\bar z})\in
SU(2)$ such that any given solution $S(t,z,{\bar z})$ to the $1+2$
dimensional Landau-Lifshitz equation (\ref{5}) is transformed to a
solution $(p(t,z,{\bar z}),q(t,z,{\bar z}),r(t,z,{\bar
z}),u(t,z,{\bar z}))$ to the following nonlinear
Schr\"odinger-type system:
\begin{eqnarray}
\left\{\begin{array}{c} iq_t-q_{z{\bar z}}+2u
q-2({\bar p}q)_z+2pq_{\bar z}+4|p|^2q=0\\
ir_t+r_{z{\bar z}}-2u r-2({\bar p}r)_z+2pr_{\bar z}-4|p|^2r=0\\
ip_t=(qr)_{\bar z}-u_z\qquad\qquad\qquad\qquad\qquad\qquad
\end{array} \right.    \label{8}
\end{eqnarray}
by the gauge transformation, where  $u$ is a unknown real function
and $p,q,r$ are unknown complex functions satisfying the following
additional restrictions
\begin{eqnarray}
{\bar p}_z+p_{\bar z}= |q|^2-|r|^2,\quad {\bar r}_z+q_{\bar
z}=2(p{\bar r}-{\bar p}q). \label{q=r}
\end{eqnarray}
\end{Theorem}

\nd{\bf Proof}. Let $S=S(t,z,{\bar z})$ be a solution to
Eq.(\ref{5}). We come to choose an $SU(2)$ matrix $G(t,z,{\bar
z})$ such that
\begin{eqnarray}
\sigma_3=-G^{-1}SG,\quad G^{-1}G_{z}=-\left(\begin{array}{cc}
  p & q\\
  {r}& -p
  \end{array}\right):=-U
  \label{14}
\end{eqnarray}
for some complex functions $p=p(t,z,{\bar z})$ and  $q=q(t,z,{\bar
z})$ and $r=r(t,z,{\bar z})$, where
$\sigma_3=\left(\begin{array}{cc}
       1& 0\\
       0& -1
       \end{array}\right)$. Indeed, by a direct computation, we
see that the general solutions to $\sigma_3=-G^{-1}SG$ are of the
forms:
\begin{eqnarray}
G=i{1\over{\sqrt{2(1-s_3)}}}(S-\sigma_3){\rm
diag}(\gamma,{\bar\gamma}), \label{gamma}
\end{eqnarray}
where $\gamma$  is a complex function of $z$, $\bar z$ and $t$ (or
in other words, $x$, $y$ and $t$) with $|\gamma|=1$.  For a fixed
$SU(2)$-matrix $G$ given in (\ref{gamma}), we have
\begin{eqnarray}
G_x=-G\left(\begin{array}{cc}
  is & \psi\\
  -{\bar \psi}& -is
  \end{array}\right), \quad G_y=-G\left(\begin{array}{cc}
  il & \phi\\
  -{\bar \phi}& -il
  \end{array}\right) \label{u}
\end{eqnarray}
for some real functions $s,l$ and complex functions $\phi, \psi$.
It is easy to see from (\ref{u}) that, for the complex variables
$z=(x+iy)/2$ and ${\bar z}=(x-iy)/2$, we have
\begin{eqnarray}
G_z=-G\left(\begin{array}{cc}
       p & (\psi-i \phi)\\
       (-{\bar \psi}+i{\bar \phi})& -p
       \end{array}\right), \quad
G_{\bar z}=-G\left(\begin{array}{cc}
       -{\bar p} & \psi+i \phi\\
       -({\bar \psi}+i{\bar \phi})& {\bar p}
       \end{array}\right),       \nn
\end{eqnarray}
where $p=l+is$ is a complex function of $x,y,t$.  Thus for any
fixed $G$ given in (\ref{gamma}) we have not only (\ref{14}) with
complex functions $p$, $q=\psi-i\phi$ and $r=-{\bar \psi}+i{\bar
\phi}$, but also
\begin{eqnarray}
G_{\bar z}=GP:=G \left(\begin{array}{cc}
       {\bar p} & {\bar r}\\
       {\bar q}& -{\bar p}
       \end{array}\right). \label{GP}
\end{eqnarray}
Furthermore,  from the integrability condition $P_{z}+U_{\bar
z}+[P,U]=0$ of the linear system: $G_{z}=-GU, G_{\bar z}=GP$, we
have
\begin{eqnarray}
{\bar p}_z+p_{\bar z}= |q|^2-|r|^2,\quad {\bar r}_z+q_{\bar
z}=2(p{\bar r}-{\bar p}q), \nn
\end{eqnarray}
which is exactly the restrictions (\ref{q=r}). Meanwhile, note
that $\al={\bar p}$ by the definition of the connection $A$ given
in (\ref{6}).

Now, for the connection $A$ given in (\ref{6}) with $S$ being
fixed above,  we define a connection 1-form as follows
\begin{eqnarray}
A^{G}&=&-G^{-1}dG+G^{-1}A G.
   \label{AG1}
\end{eqnarray}
From $S=-G\sigma_3G^{-1}$, $G_{\bar z}=GP$ and $\al={\bar p}$, we
have
\begin{eqnarray}
&&SS_{\bar z}+2{\al} S\nn \\
&=& -G\sigma_3G^{-1}(-GP\sigma_3G^{-1}+G
\sigma_3PG^{-1})+2{\bar p}(-G\sigma_3G^{-1})\nn\\
&=&-2GP^{(\hbox{off-diag})}G^{-1}-2{\bar p}G\sigma_3G^{-1}\nn\\
&=& -2G_{\bar z}G^{-1}.  \label{GG}
\end{eqnarray}
So (\ref{AG1}) can be re-expressed as follows
\begin{eqnarray}
A^{G} &=& (-G^{-1}G_{\bar z}+G^{-1}VG)d{\bar
z}+(i\lambda\sigma_3+U)d{z}\nn\\
&&+\bigg(-G^{-1}G_t+ \lambda(2G^{-1}VG-2G^{-1}G_{\bar z})\bigg)dt.
   \label{AG}
\end{eqnarray}
Since $A$ satisfies the prescribed curvature condition:
$$
F_A=dA-A\wedge A=K,
$$
where $K$ is given by (\ref{K}), from gauge theory we know that
$A^G$ must satisfies
\begin{eqnarray}
F_{A^{G}}=dA^{G}-A^{G}\wedge A^{G}=G^{-1}KG. \label{017}
\end{eqnarray}
Comparing respectively the coefficients of $d{z}\wedge d{\bar z}$,
$d{\bar z}\wedge dt$ and $d{z}\wedge dt$ in the both sides of
(\ref{017}), we have
\begin{eqnarray}
\bigg(-G^{-1}G_{\bar z}+G^{-1}VG\bigg)_{z}-U_{\bar
z}-\bigg[i\lambda \sigma_3+U, -G^{-1}G_{\bar z}+G^{-1}VG\bigg]=0,
\label{FA1}
\end{eqnarray}
\begin{eqnarray}
&&-\bigg(-G^{-1}G_{\bar z}+G^{-1}VG\bigg)_t+\bigg(-G^{-1}G_t+
\lambda(2G^{-1}VG-2G^{-1}G_{\bar z})\bigg)_{\bar z}
\nn \\
&&-\bigg[-G^{-1}G_{\bar z}+G^{-1}VG,-G^{-1}G_t+
\lambda(2G^{-1}VG-2G^{-1}G_{\bar z})\bigg]
\nn\\
&=&-G^{-1}V_tG+2\lambda G^{-1}V_{\bar z}G+\lambda
G^{-1}\bigg(SS_{\bar z}+2{\al}
S\bigg)_{\bar z}G\nn\\
&&-G^{-1}\bigg[V,2\lambda V+\lambda(SS_{\bar z}+2{\al} S)\bigg]G,
\label{FA2}
\end{eqnarray}
and
\begin{eqnarray}
&&-U_t+\bigg(-G^{-1}G_t+ \lambda(2G^{-1}VG-2G^{-1}G_{\bar
z})\bigg)_{z}\nn
\\&& -\bigg[i\lambda\sigma_3+U,-G^{-1}G_t+
\lambda(2G^{-1}VG-2G^{-1}G_{\bar z})\bigg]\nn\\
&=&G^{-1}\lambda\bigg({1\over2}[S_{z},S_{\bar z}]+2(\al
S)_z\bigg)G. \label{FA3}
\end{eqnarray}
Setting
\begin{eqnarray}{\wt V}=-G^{-1}G_{\bar z}+G^{-1}VG,  \label{wV}
\end{eqnarray}
and by a straightforward computation, we obtain
\begin{eqnarray}
 G^{-1}\bigg((i\lambda S)_{\bar z}+V_{z} +[i\lambda S,V]
\bigg)G=
 -U_{\bar z}+ {\widetilde V}_{z}-[i\lambda \sigma_3+U,{\widetilde
V}],  \nn %\label{UV}
\end{eqnarray}
which implies that (\ref{FA1}) is automatically satisfied from
(\ref{V}). By using the definition ${\wt V}=-G^{-1}G_{\bar
z}+G^{-1}VG$ again, we see that (\ref{FA2}) is equivalent to
\begin{eqnarray}
&&-{\widetilde V}_t+\bigg(-G^{-1}G_t+2\lambda {\widetilde
V}\bigg)_{\bar z} -[{\widetilde V}, -G^{-1}G_t+2\lambda
{\widetilde
V}]\nn\\
&=&G^{-1}\bigg\{-V_t+\lambda (2V +SS_{\bar z}+2{\al} S)_{\bar
z}-[V,\lambda(2V +SS_{\bar z}+2{\al} S)]\bigg\}G. \label{FA4}
\end{eqnarray}
It is also a straightforward verification that (\ref{FA4}) is an
identity too. Finally, we come to treat (\ref{FA3}). First, by
using  the identities: $S_z=G(U\sigma_3-\sigma_3 U)G^{-1}$ and
$S_{\bar z}=-GP\sigma_3 G^{-1}+G\sigma_3PG^{-1}$ deduced from
(\ref{14}) and (\ref{GP}) respectively, $\al={\bar p}$ and the
first equation of (\ref{q=r}), we have
\begin{eqnarray}
&&{{1\over2}[S_{z},S_{\bar z}]+2(\al S)_{z}} \nn\\
&&=G2\left([U^{\hbox{(off-diag)}}\sigma_3,\sigma_3P^{\hbox{(off-diag)}}]-2
{\bar p}_{z}\sigma_3+4{\bar p} U^{\hbox{(off-diag)}}\sigma_3\right)G^{-1}\nn\\
&&=G\bigg(2p_{\bar z}\sigma_3+4{\bar p}
U^{\hbox{(off-diag)}}\sigma_3\bigg)G^{-1}=G\bigg(2p_{\bar
z}\sigma_3+ [\sigma_3,{Q}]\bigg)G^{-1}, \label{dxt}
\end{eqnarray}
%since $[\sigma_3,Q]=4{\bar p} U^{\hbox{(off-diag)}}\sigma_3$,
where ${ Q}=\left(\begin{array}{cc}
       0 & -2q{\bar p}\\
       -2{r}{\bar p}& 0
       \end{array}\right)$.
Thus, (\ref{FA3}) is equivalent to
\begin{eqnarray}
-U_t+(-G^{-1}G_t+2\lambda {\wt V})_{z}
-[i\lambda\sigma_3+U,-G^{-1}G_t+ \lambda2{\wt V}] =2\lambda
p_{\bar z}\sigma_3 +\lambda[\sigma_3,{ Q}], \nn
\end{eqnarray}
or equivalently,
\begin{eqnarray}
&&-U_t+(-G^{-1}G_t)_{z}+[U,G^{-1}G_t] \nn\\ && +
\lambda\left(2U_{\bar z}-2p_{\bar z}\sigma_3+
i[\sigma_3,G^{-1}G_t]-[\sigma_3,{ Q}]\right)=0 \label{EQ}
\end{eqnarray}
here we have used the identities (\ref{FA1}): $U_{\bar z}-{\wt
V}_{z}+[i\lambda \sigma_3+U, {\wt V}]=0$. Comparing the
coefficients of $\lambda$ and the constant term in (\ref{EQ}), we
obtain
\begin{eqnarray}
&&-U_t+(-G^{-1}G_t)_{z}+[U,G^{-1}G_t]=0,\label{term0} \\
&&2U_{\bar z}-2p_{\bar
z}\sigma_3+i[\sigma_3,G^{-1}G_t]-[\sigma_3,{ Q}]=0. \label{term1}
\end{eqnarray}
The vanishing of the diagonal part of (\ref{term0}) and the
equation (\ref{term1}) lead to
\begin{eqnarray}
&&G^{-1}G_t=-i\left\{\left(u+U_{\bar
z}^{(\hbox{off-diag})}\right)\sigma_3+Q\right\}  \label{019} \\
&&p_t=i(u_z-(qr)_{\bar z})
%=i((qr)_z-(qr)_{\bar z})-i\tau_z
\label{020}
\end{eqnarray}
for some real function $u=u(t,x,y)$. Here we have used the second
equation of (\ref{q=r}) to verify that $-i(U_{\bar
z}^{(\hbox{off-diag})}\sigma_3+Q)\in su(2)$. Substituting
(\ref{wV}) and (\ref{019},\ref{020}) into (\ref{AG}) and
(\ref{017}) respectively, we obtain
\begin{eqnarray}
{A^G}= {\widetilde V}d{\bar z}+\bigg(i\lambda\sigma_3+U\bigg)d{z}
+ \left(2\lambda {\widetilde V}+i\left(u +U_{\bar
z}^{(\hbox{off-diag})}\right)\sigma_3+iQ\right) dt
\end{eqnarray}
and
\begin{eqnarray}
{K}^G=G^{-1}KG&=& \left\{-{\widetilde V}_t+\left(2\lambda
{\widetilde V}+iu\sigma_3 +iU_{\bar
z}^{(\hbox{off-diag})}\sigma_3+iQ\right)_{\bar z}\right.\nn\\
&&\left.-[{\widetilde V}, 2\lambda {\widetilde V}+iu\sigma_3
+iU_{\bar z}^{(\hbox{off-diag})}\sigma_3+iQ]\right\}d{\bar
z}\wedge
dt\nn\\
&& +\lambda\bigg(2U^{(\hbox{diag})}_{\bar
z}+[\sigma_3,Q]\bigg)d{z}\wedge dt
\end{eqnarray}
where $\lambda$ is the same spectral parameter as in (\ref{6}),
${\widetilde V}={\widetilde V}(\lambda, z,\bar z,t)$ is a
$2\times2$-matrix satisfying
\begin{eqnarray}
U_{\bar z}- {\widetilde V}_{z}+[i\lambda \sigma_3+U,{\widetilde
V}]=0. \label{tildeV}
\end{eqnarray}
So, it is a direct computation, from the prescribed curvature
representation: $F_{A^G}=K^G$, that the corresponding equation for
unknown functions $(p,q,r,u)$ is just the nonlinear
Schr\"odinger-type system (\ref{8}). We would like to point out
that one should apply the identity (\ref{tildeV}) in the
computation. This completes the proof of Theorem 1.$\Box$

\bigskip
Now we come to consider the nonlinear Schr\"odinger-type system
(\ref{8})  with the restriction (\ref{q=r}). As indicated in
Theorem 1, it is a PDE with prescribed curvature representation:
\begin{eqnarray}
F_{\widetilde A}=d{\widetilde A}-{\widetilde A}\wedge {\widetilde
A} ={\widetilde K} \label{09}
\end{eqnarray}
in which
\begin{eqnarray}
{\widetilde A}:= {\widetilde V}d{\bar
z}+\bigg(i\lambda\sigma_3+U\bigg)d{z} + \bigg\{2\lambda
{\widetilde V}+i(u +U_{\bar z}^{(\hbox{off-diag})})\sigma_3+i
Q\bigg\} dt \label{9}
\end{eqnarray}
and
\begin{eqnarray}
{\widetilde K}&:=& \Bigg\{-{\widetilde V}_t+\left(2\lambda
{\widetilde V}+i(u
+U_{\bar z}^{(\hbox{off-diag})})\sigma_3+i Q\right)_{\bar z}\nonumber\\
&&-[{\widetilde V}, 2\lambda {\widetilde V}+i(u +U_{\bar
z}^{(\hbox{off-diag})})\sigma_3+i Q]\Bigg\}d{\bar z}\wedge dt
\nn\\
&& +\lambda\bigg(2U^{(\hbox{diag})}_{\bar
z}+[\sigma_3,Q]\bigg)d{z}\wedge dt, \label{TK}
\end{eqnarray}
where all the notations have the same meanings indicated above,
i.e., $U=\left(\begin{array}{cc}
  p & q\\
  {r}& -p
  \end{array}\right)$, $\widetilde V$ solves the equation
  (\ref{tildeV}) and $Q=\left(\begin{array}{cc}
       0 & -2q{\bar p}\\
       -2{r}{\bar p}& 0
       \end{array}\right)$.

Next we shall prove that the above gauge transformation from the
$1+2$ dimensional isotropic Landau-Lifshitz (\ref{5}) to the
nonlinear Schr\"odinger-type system (\ref{8}) with the
restrictions (\ref{q=r}) is in fact reversible.

\begin{Theorem}
There is a gauge matrix function $G(t,x,y)\in SU(2)$ such that any
solution $(p,q,r,u)$ to the nonlinear Schr\"odinger-type system
(\ref{8}) with the restriction (\ref{q=r}) can be transformed to a
solution $S$ to the $1+2$ dimensional isotropic Landau-Lifshitz
equation (\ref{5}) by the gauge transformation of $G$. Moreover,
if we require that the gauge matrix $G$ satisfies
$G|_{x={y}=t=0}=I$. Then any $C^m$-solution ($m\ge 2$) to the
nonlinear Schr\"odinger-type system (\ref{8}) with the restriction
(\ref{q=r}) corresponds uniquely to a $C^{m+1}$-solution to the
Schr\"odinger map (\ref{5}) and vice versa.
\end{Theorem}

\noindent{\bf Proof}:  Let $(p,q,r,u)$ be a solution to
Eq.(\ref{8}) with the restrictions (\ref{q=r}). From the
prescribed curvature formulation (\ref{09}), it is a key
observation that Eq.(\ref{8}) with the restrictions (\ref{q=r}) is
in fact the integrability condition of the following linear
system:
\begin{eqnarray}
G_{z}=-GU,\quad G_t=-Gi\left((u +U^{(\hbox{off-diag})}_{\bar
z})\sigma_3+Q\right),\label{10}
\end{eqnarray}
or equivalently,
\begin{eqnarray}
\left\{\begin{array}{c}
G_{x}=-G\left(\begin{array}{cc}
  i\hbox{Im}{p} & \psi\\
  -{\bar \psi}& -i\hbox{Im}{p}
  \end{array}\right)\qquad\\
 G_y=-G\left(\begin{array}{cc}
  i\hbox{Re}{p} & \phi\\
  -{\bar \phi}& -i\hbox{Re}{p}
  \end{array}\right)\qquad\\
 G_t=G\left(\begin{array}{cc}
       -iu  & iq_{\bar z}+2i{\bar p}q\\
       -i{r}_{\bar z}+2i{\bar p}{r}& iu
       \end{array}\right)
       \end{array}\right.
       \label{r10}
\end{eqnarray}
where $\psi=(q-{\bar r})/2$ and $\phi=i(q+{\bar r})/2$ (here we
have used (\ref{q=r}) to verify that the coefficient matrix in
righthand side of the third equation of (\ref{r10}) is an
$su(2)$-matrix). This implies that (\ref{10}) or (\ref{r10}) is a
compatible linear differential system. Since the coefficient
matrices in (\ref{r10}) are $su(2)$-matrices, it indicates that
general solutions $G(t,x,y)$ to (\ref{r10}) or equivalently
(\ref{10}) belong to $SU(2)$ group. Now let $G(t,z,{\bar
z})=G(t,x,y)\in SU(2)$ be a fundamental solution to (\ref{10}) or
equivalently (\ref{r10}), and we consider the following gauge
transformation,
\begin{eqnarray}
A=(dG) G^{-1} +G {\widetilde A}G^{-1}, \label{11}
\end{eqnarray}
where ${\widetilde A}$ is the 1-form connection given in (\ref{9})
with $(p,q,r,u)$ being given above. We try to show that the 1-form
$A$ defined by (\ref{11}) is exactly the connection of
Eq.(\ref{5}) given in (\ref{6}) when $S$ and $\al$ are suitably
determined. In fact, substituting the coefficient $-i\lambda S$ of
$d{z}$ of (\ref{6}) into (\ref{11}) and comparing the coefficients
of $\lambda$ of $d{z}$ in the both sides of (\ref{11}), we obtain
\begin{eqnarray}
G_{z}=-GU,\quad S=-G\sigma_3G^{-1} \quad ({\rm hence} \quad
S^2=I). \label{12}
\end{eqnarray}
The first equation of (\ref{12}) is automatically satisfied
because of the first equation of (\ref{10}). The second one of
(\ref{12}) is regarded as defining $S$. Now, we have to prove that
the coefficients of $d{\bar z}$ and $dt$ of $A$ defined by
(\ref{11}) are respectively the same coefficients of $d{\bar z}$
and $dt$ of the connection given in (\ref{6}), that is,
\begin{eqnarray}
 V&=&G_{\bar z}G^{-1}+G{\wt V} G^{-1},  \label{dy} \\
\lambda\bigg(2V +SS_{\bar z}+2\al S\bigg)
&=&G_tG^{-1}+G\left(2\lambda {\widetilde
V}+i\bigg((u+U^{(\hbox{off-diag})}_{\bar
z})\sigma_3+Q\bigg)\right)G^{-1}. \label{dt}
\end{eqnarray}
Eq.(\ref{dy}) can be regarded as defining $V$ if we can show that
such a $V$  solves Eq.(\ref{V}), i.e., for the $S$ being given in
(\ref{12}) we have
\begin{eqnarray}
 (-i\lambda S)_{\bar z}-V_{z} +[-i\lambda S,V]=0. \label{SV}
\end{eqnarray}
The proof of (\ref{SV}) is a direct computation. Indeed, by using
the expression of $V$ given in (\ref{dy}) and the fact that $G$
fulfills (\ref{r10}) (this equivalent to having (\ref{14}),
(\ref{GP}) and the third equation of (\ref{r10})), we have
\begin{eqnarray}
 (i\lambda S)_{\bar z}+V_{z} +[i\lambda S,V]=
 G\left(U_{\bar z}- {\widetilde V}_{z}+[i\lambda \sigma_3+U,{\widetilde
V}] \right)G^{-1}. \label{dy=dy}
\end{eqnarray}
Since ${\wt V}$ satisfies (\ref{tildeV}), this establishes
(\ref{SV}).  For proving (\ref{dt}), since $G$ satisfies the
second equation of (\ref{10}), it is easy to see that the proof of
(\ref{dt}) is equivalent to
\begin{eqnarray}
2V +SS_{\bar z}+2{\al} S=2G{\wt V}G^{-1},  \nn
\end{eqnarray}
which, because of (\ref{dy}), is  equivalent to
\begin{eqnarray}
SS_{\bar z}+2{\al} S=-2G_{\bar z}G^{-1}.  \label{dt=dt}
\end{eqnarray}
Now we take $\al={\bar p}$ which fulfills the requirement of $\al$
in the definition of the connection (\ref{6}). From $G_{z}=-GU$
$G_{\bar z}=GP$ and $\al={\bar p}$,  we have
\begin{eqnarray}
&&-(SS_{\bar z}+2{\al} S)\nn \\
&=& -G\sigma_3G^{-1}(-GP\sigma_3G^{-1}+G
\sigma_3PG^{-1})+2(-{\bar p})(-G\sigma_3G^{-1})\nn\\
&=&2GP^{(\hbox{off-diag})}G^{-1}+2{\bar p} G\sigma_3G^{-1}\nn\\
&=& 2GPG^{-1}= 2G_{\bar z}G^{-1}. \nn% \label{G}
\end{eqnarray}
This proves (\ref{dt=dt}). Thus we have proved that the two
connections given by (\ref{11}) and (\ref{6}) respectively are
actually the same one when $S=-G\sigma_3G^{-1}$ and $\al={\bar
p}$. What's the remainder for us to do is to prove that the
curvature formula
\begin{eqnarray}
F_A=K =G{\widetilde K}G^{-1}=GF_{\widetilde A}G^{-1} \label{curv}
\end{eqnarray}
under the gauge transformation is satisfied too, where $K$ is
given by (\ref{K}) and ${\widetilde K}$ is given by (\ref{TK}). In
fact, on the one hand, we see that
\begin{eqnarray}
G{\wt K}G^{-1}&=& G\Bigg(\bigg\{-{\widetilde V}_t+\left(2\lambda
{\widetilde V}+i(u+U^{(\hbox{off-diag})}_{\bar z})
\sigma_3+iQ\right)_{\bar z}\nonumber\\
&&-[{\widetilde V}, 2\lambda {\widetilde
V}+i(u+U^{(\hbox{off-diag})}_{\bar
z})\sigma_3+i Q]\bigg\}d{\bar z}\wedge dt\nn\\
&&+\lambda(2U^{(\hbox{diag})}_{\bar z}+[\sigma_3,Q])d{z}\wedge
dt\Bigg)G^{-1}. \label{GKG}
\end{eqnarray}
On the other hand, by using (\ref{dy}, \ref{dt=dt}) and
(\ref{10}), it is a straightforward calculation that the
coefficient of $d{\bar z}\wedge dt$ in $K$ is
\begin{eqnarray}
&&-V_t+\lambda (2V +SS_{\bar z}+2{\al} S)_{\bar z}-[V,\lambda(2V
+SS_{\bar z}+2{\al} S)]\nonumber \\
&=&-V_t+2\lambda(G{\wt V}G^{-1})_{\bar z}-[V,2\lambda G{\wt
V}G^{-1}]\nn \\
&=&-G_{t{\bar z}}G^{-1}+G_{\bar z}GG_tG^{-1}-G_t{\wt
V}G^{-1}-G{\wt
V}_tG^{-1}+G{\wt V}G^{-1}G_tG^{-1}+2\lambda G_{\bar z}{\wt V}G^{-1}\nn\\
&&+2\lambda G{\wt V}_{\bar z}G^{-1}-2\lambda G{\wt V}G^{-1}G_{\bar
z}G^{-1}-2\lambda(G_{\bar z}{\wt V}G^{-1}-G{\wt
V}G^{-1}G_{\bar z}G^{-1})\nn\\
&=&G\Bigg\{-{\widetilde V}_t+\left(2\lambda {\widetilde
V}+i(u+U^{(\hbox{off-diag})}_{\bar z})\sigma_3
+iQ\right)_{\bar z}\nonumber\\
&&-[{\widetilde V}, 2\lambda {\widetilde
V}+i(u+U^{(\hbox{off-diag})}_{\bar z})\sigma_3+i Q]\Bigg\}G^{-1}.
\label{dydt}
\end{eqnarray}
(\ref{GKG}) and (\ref{dydt}) indicate that the two coefficients of
$d{\bar z}\wedge dt$ in the both sides of (\ref{curv}) are the
same one. Meanwhile, by applying the similar argument in getting
(\ref{dxt}), we have
\begin{eqnarray}
{1\over2}[S_{z},S_{\bar z}]+2(\al S)_{z}=G(2p_{\bar z}\sigma_3+
[\sigma_3,{Q}])G^{-1}, \label{dxdt}
\end{eqnarray}
where ${ Q}=\left(\begin{array}{cc}
       0 & -2q{\bar p}\\
       -2{r}{\bar p}& 0
       \end{array}\right)$ as before. (\ref{dxdt}) implies
\begin{eqnarray}
\lambda\bigg({1\over2}[S_{z},S_{\bar z}]+2(\al
S)_{z}\bigg)=\lambda G \bigg(2U^{(\hbox{diag})}_{\bar
z}+[\sigma_3,Q]\bigg)G^{-1}, \nn
\end{eqnarray}
which shows that the two coefficients of $d{z}\wedge dt$ in $K$
and $G{\wt K}G^{-1}$ are also the same one. Thus we have proved
the desired identity (\ref{curv}), which implies the holding of
the prescribed curvature representation (\ref{7}). Hence we obtain
that $S$ solves Eq.(\ref{5}). This indicates that $S$ defined by
the second equation of (\ref{12}) from a solution $(p,q,r,u)$ to
(\ref{8}) with the restrictions (\ref{q=r}) satisfies the 1+2
dimensional Landau-Lifshitz equation (\ref{5}).

Since $G$ is a solution to the linear first-order differential
system (\ref{10}), it is well-known from linear theory of
differential equations that such a $G$ is unique if we propose the
initial condition $G|_{x=y=t=0}=I$ on $G$. Under this
circumstance, we see that a solution $(p,q,r,u)$ to Eq.(\ref{8})
with the restrictions (\ref{q=r}) corresponds uniquely to a
solution $S$ to Eq.(\ref{5}) by the gauge transformation and vice
versa. Furthermore, because of the relation
$S_z=G2U^{(\hbox{off-diag})}\sigma_3G^{-1}$ deduced from
(\ref{12}), the remainder part of the theorem is obviously true.
$\Box$

\bigskip

%\begin{Remark}
We would like to point out that the unknown
functions in the system (\ref{8}) can be reduced to $(p,q,r)$ with
$p$ being a real function. In fact, we may restrict $G$ given by
(\ref{gamma}) to satisfy
$
G_x=-G\left(\begin{array}{cc}
  0 & \psi\\
  -{\bar \psi}& 0
  \end{array}\right)$ for some complex function $\psi$ which leads to
$\gamma=\exp\left(-i{\int^{x}_0{{(s_1{s_2}_{x}
-s_2{s_1}_{x})}\over{2s_3-2}}d{x} }\right)$.  Thus  $p$ becomes
now a real function of $x,y$ and $t$. Under this situation, it is
easy to verify that the vanishing of the diagonal part of
(\ref{term0}) and the equation (\ref{term1}) lead to
\begin{eqnarray}
&&G^{-1}G_t=-i\left\{\left((\hbox{Re}(qr)-\pa^{-1}_x\pa_y\hbox{Im}(qr))+U_{\bar
z}^{(\hbox{off-diag})}\right)\sigma_3+Q\right\}+i\tau\sigma_3  \label{0019} \\
&&p_t=2\bigg(\hbox{Re}(qr)\bigg)_y+\left((\hbox{Im}(qr))_x
-\pa^{-1}_x\pa_y^2\hbox{Im}(qr)\right)- \tau_y \label{0020}
\end{eqnarray}
for some real function $\tau=\tau(y,t)$ depending only on $t$ and
$y$. Moreover, notice that the restriction (\ref{14}) on the gauge
matrix $G$ allows an arbitrariness in $G$ of the form: $G\to G
e^{i\sigma_3\beta(t,y)}$ for an arbitrary real function
$\beta(t,y)$. If we require $\beta$ (the existence of such a
$\beta$ is easy to verify) to satisfy
\begin{eqnarray}
{{\partial\beta}\over {\partial y}}=-\int^t_0\tau_y, \quad
{{\partial\beta}\over {\partial t}}=-\tau \nn % \label{ga2}
\end{eqnarray}
then $G$ can be modified so that for the new $G$ we have
\begin{eqnarray}
G_{y}&=&-G\left(\begin{array}{cc}
       i p+i\int^t_0\tau_y & \psi e^{-2i\beta}\\
       -{\bar \psi}e^{2i\beta}& -ip-i\int^t_0\tau_y
       \end{array}\right), \nn \\
G_{t}&=&G\left(\begin{array}{cc}
       -i({\hbox{Re}}(qr)-\partial^{-1}_x
\partial_y({\hbox{Im}}(qr))) -i\tau & \sigma e^{-2i\beta}\\
       -{\bar \sigma}e^{2i\beta}& i({\hbox{Re}}(qr)-\partial^{-1}_x
\partial_y({\hbox{Im}}(qr)))+i\tau
       \end{array}\right), \nn
\end{eqnarray}
where $\sigma=iq_{\bar z}+2i{p}q$, which implies that for the new
$G$ the second term on the right-hand side of (\ref{0019}) is $0$
and meanwhile the third term on the right-hand side of
(\ref{0020}) is $0$ too. Hence the system (\ref{8}) is reduced to
the following nonlinear Schr\"odinger type system:
\begin{eqnarray}
\left\{\begin{array}{c} iq_t-q_{z{\bar
z}}+2q(({\hbox{Re}}(qr)-\partial^{-1}_x
\partial_y({\hbox{Im}}(qr)))+2pq_{\bar z}-2(pq)_z+4p^2q=0\\
ir_t+r_{z{\bar z}}-2r(({\hbox{Re}}(qr)-\partial^{-1}_x
\partial_y({\hbox{Im}}(qr)))+2pr_{\bar z}-2(pr)_z-4p^2r=0\\
p_t=2({\hbox{Re}}(qr))_y+({\hbox{Im}}(qr))_x-\partial^{-1}_x
\partial_y^2({\hbox{Im}}(qr)),\quad\qquad\qquad\qquad\qquad
\end{array} \right.    \label{r28}
\end{eqnarray}
where the real function $p$ and the complex function $q,r$ must
satisfy the  restriction (\ref{q=r}) too. Though (\ref{r28}) looks
much simpler than (\ref{8}) in unknown variables, both (\ref{8})
and (\ref{r28}) with the restrictions (\ref{q=r}) are essentially
equivalent to each other. However, the natural choice of complex
version of $p$ in the system (\ref{8}) plays an important role in
constructing blow-up $H^3({\bf R}^2)$-solutions to the
Landau-Lifshitz equation (\ref{4}), as we shall see in the next
section.
%\end{Remark}

\begin{Remark}
The expression (\ref{14}) or (\ref{12}) of the two gauge
equivalent solutions gives an explicit relationship of the 1+2
isotropic dimensional Landau-Lifshitz equation (\ref{4}) and the
nonlinear Schr\"odinger-type system (\ref{8}), which also plays a
key role in constructing blow-up solutions to Eq.(\ref{4}) in the
next section. This explicit relation is not obtained by other
transformations, such as the Hasimoto transformation.
\end{Remark}

\section * {\S 3. Blowing up solutions}

We follow the basic conventional notations for Sobolev spaces
$H^{k}({\bf R}^n)$, $W^{k,\sigma}({\bf R}^n)$ ($\sigma>1$) of real
or complex-valued functions or spaces $C^k({\bf R}^n)$ of
continuous differential functions up to order $k$ on ${\bf R}^n$
for $n\ge2$ and norms $||\cdot||_{W^{k,\sigma}({\bf R}^n)}$ or
$||\cdot||_{C^k({\bf R}^n)}$ used in \cite{GiT}. In this section,
we shall construct, by use of its gauge equivalent nonlinear
Schr\"odinger-type equation (\ref{8}) displayed in the previous
section, blowing up $H^3({\bf R}^2)$-solutions to the
Landau-Lifshitz equation (\ref{4}). Before doing this, let's
characterize some anayltic properties of the system (\ref{8}) or
equivalently Eq.(\ref{r28}) with the restrictions (\ref{q=r}).

\bigskip

\nd{\bf Claim 1}. System (\ref{r28}) has the following
conservation laws:
\begin{eqnarray}
&&\int_{{\bf R}^2}|q(x,y,t)|^2dxdy=\int_{{\bf
R}^2}|q(x,y,0)|^2dxdy,\nn\\
&&\int_{{\bf R}^2}|r(x,y,t)|^2dxdy=\int_{{\bf
R}^2}|r(x,y,0)|^2dxdy.\nn
\end{eqnarray}
In fact, we multiply $2{\bar q}$ (or $2{\bar r}$) to the first
equation (or the second equation) of (\ref{r28}) and take the
imaginary part of the result to get
$$
{{\pa}\over{\pa t}}|q|^2=2\hbox{Im}((\Delta q){\bar q})
-4p(|q|^2)_y-4p_y|q|^2 \quad ({{\pa}\over{\pa
t}}|r|^2=-2\hbox{Im}((\Delta r){\bar r})-4p(|r|^2)_y-4p_y|r|^2).
$$
Thus the above conservation laws follow from integrating the
identities over ${\bf R}^2$.

\bigskip
\nd{\bf Claim 2}. Let us introduce the polar coordinates
$(\rho,\theta)$ of ${\bf R}^2$, that is $x=\rho\cos\theta,
y=\rho\sin\theta$. Thus we have
$${{\partial}\over{\partial z}}=
{{e^{-i\theta}}}\left({{\partial}\over{\partial \rho}}-i{1\over
\rho}{{\partial}\over{\partial \theta}}\right),\quad
{{\partial}\over{\partial {\bar z}}}=
{{e^{i\theta}}}\left({{\partial}\over{\partial \rho}}+i{1\over
\rho}{{\partial}\over{\partial \theta}}\right).$$ We would like to
find following ansatz solutions to (\ref{8}):
\begin{eqnarray}
p=0, q=e^{-i\theta}Q(\rho,t), r=-e^{-i\theta}{\bar Q}(\rho,t),
u=-|Q|^2+2\int_{\rho}^{\infty}{{|Q|^2(\tau,t)}\over{\tau}}d\tau
\end{eqnarray}
for some suitable functions $Q(\rho,t)$. One may verify that
(\ref{q=r}) and the third equation of (\ref{8}) are satisfied
automatically (since $
\partial^{-1}_z\partial_{\bar z}(rq)=\partial^{-1}_z\partial_{\bar
z}(-e^{-2i\theta}|Q|^2(\rho,t))=-|Q|^2(\rho,t)+2\int_{\rho}^{\infty}
{{|Q|^2(\tau,t)}\over{\tau}}d\tau $) and the first and second
equations of (\ref{8}) lead to
\begin{eqnarray}
 iQ_t-\left(Q_{\rho\rho}+{1\over{\rho}}Q_{\rho}
-{1\over{\rho^2}}Q\right)-2Q\left(|Q|^2 -2\int_{\rho}^{\infty}
{{|Q(\tau,t)|^2}\over{\tau}}d\tau \right)=0. \label{QQ}
\end{eqnarray}
This equation, or its equivalent form:
$iQ_t+Q_{\rho\rho}+{1\over{\rho}}Q_{\rho}
-{1\over{\rho^2}}Q+2Q\left(|Q|^2 -2\int_{\rho}^{\infty}
{{|Q(\tau,t)|^2}\over{\tau}}d\tau \right)$ $=0$, was deduced by
many authors from the (generalized) Hasimoto transformation (see
\cite{La,bib:LD,bib:RJ, CSU}).

\bigskip
\nd{\bf Claim 3}. If $q(t,s)$ is an arbitrary solution to the
following nonlinear Schr\"odinger equation:
\begin{eqnarray}
iq_t-q_{ss}-2q|q|^2=0
    \label{nn}
\end{eqnarray}
which is an integrable system, one can verify straightforwardly
that $(p=0,q(t,\cos\delta z+\sin\delta{\bar z}),$ $r(t,\cos\delta
z+\sin\delta{\bar z})=-{{1-i}\over{1+i}}{\bar q}, u=-|q|^2)$ is a
solution to (\ref{8}) with the restriction (\ref{q=r}), where
$\delta$ is a free parameter. It is well-known that system
(\ref{nn}) has (global) $1+1$ dimensional $N$-soliton solutions
(see \cite{FaT}) and hence so does (\ref{4}) correspondingly. For
example, the following 1+2 dimensional travelling 1-soliton to
Eq.(\ref{4}),
\begin{eqnarray}
&&s_1(x,y,t)={{2\sh(\cos \delta x+\sin \delta y)}\over{\ch^2(\cos
\delta
x+\sin\delta y)}}\cos t\nn \\
&&s_2(x,y,t)={{2\sh(\cos\delta x+\sin\delta
y)}\over{\ch^2(\cos\delta x
+\sin\delta y)}}\sin t\nn \\
&&s_3(x,y,t)=1-{{2}\over{\ch^2(\cos\delta x+\sin \delta y)}},\nn
\end{eqnarray}
is obtained from the 1-soliton solution $q=e^{-it}\hbox{sech} s$
to (\ref{nn}) by the gauge transformation.

\bigskip
The fact that the Landau-Lifshitz equations are related to
nonlinear Schr\"odinger-type equations has been known for a long
time. Many authors applied properties of nonlinear
Schr\"odinger-type equations to study the Landau-Lifshitz
equations. For example, it was the use of its equivalent nonlinear
Schr\"odinger-type equation obtaining by the technique of the
stereographic projection $S^2\to C$, Sulem, Sulem and Bardos
proved in \cite{SuSu} the global $W^{m+1,6}({\bf R}^n)$-existence
of the Cauchy problem of the Landau-Lifshitz eqaution (\ref{0})
(with $n\ge2$) for small initial data. So it is very reasonable
that we may use the nonlinear Schr\"odinger-type system (\ref{8})
to reveal the blow-up phenomenon of the Landau-Lifshitz equation
(\ref{4}), though this nonlinear Schr\"odinger-type system looks
very complicated. Let briefly review some blow-up results and
searching methods of nonlinear Schr\"odinger equations since they
will enlighten  on our approach. There has been much interest and
ground-breaking work within the decades in the study of nonlinear
Schr\"odinger equations with general nonlinearities (see, for
example, \cite{B1,B2,Ke,SuS}).  Blowing-up solutions to the Cauchy
problem of nonlinear Schr\"odinger equations of the forms
\begin{eqnarray}
\left\{\begin{array}{c}
iq_t=\Delta q+|q|^{\sigma-1}q, \quad {\rm in} \quad  {\bf R}^n\\
q(x,0)=q_0(x), \qquad\qquad\qquad
\end{array}
\right., \label{eq}
\end{eqnarray}
where $\sigma>1+n/4$ is a positive constant, were displayed in
\cite{Gla,Ka}. In the proofs of the mentioned blow-up results for
the nonlinear Schr\"odinger equation (\ref{eq}), the conservation
law $E(q)={1\over 2}\int_{{\bf R}^n}|\nabla
q_0|^2dv-{1\over{\sigma+1}} \int_{{\bf
R}^n}|q_0|^{\sigma+1}dv=E(q_0)$ plays a key role. Comparing the
nonlinear Schr\"odinger equation (\ref{eq}) with the nonlinear
Schr\"odinger-type equation (\ref{8}) or (\ref{QQ}), we find that
there are additional integral terms in our present case. The extra
term(s) prevent us from getting the analogous conservation law as
that of Eq.(\ref{eq}). Thus we will face new difficulties in
characterizing the blowing-up property of the present nonlinear
Schr\"odinger-type system if we go along the way depending on
conservation laws. On the other hand, the Eq.(\ref{eq}) of the
critical case $\si=4/n+1$ admits the following conformal
invariance (\cite{Wei}),
$$
q(x,t)\to
(Cq)(x,t)={{e^{-i{{b|x|^2}\over{4(a+bt)}}}}\over{(a+bt)^{n/2}}}q(X,T)
$$
where $X={{x}\over{a+bt}}$, $T={{c+dt}\over{a+bt}}$ and
$\left(\begin{array}{cc}
              a & b\\
              c& d\end{array} \right)\in SL_2({\bf R})$, i.e.
              $a, b, c, d$ are real numbers and $ad-bc=1$. That is
              to say $Cq$ is a solution to (\ref{eq}) too.
Weinstein constructed in \cite{Wei} (one can also refer to
\cite{B2}) by using this conformal invariance (which is absent for
$\si\not=4/n+1$) blowing-up solutions from localized finite energy
solitary waves. And he also proved the sharpness of a condition
for the global existence of solutions in \cite{bib:wein}. However,
we find surprisingly that the method used by Weinstein can be
modified to our present nonlinear Schr\"odinger-type equations
(\ref{8}) and applied to construct blowing-up solutions to the
Landau-Lifshitz equation (\ref{4}).

\begin{Lemma}
Assume that $q=q(\rho,t)$, where $(\rho,\theta)$ is the polar
coordinates of ${\bf R}^2$,  solves the following nonlinear
Schr\"odinger-type equation:
\begin{eqnarray}
iq_t-q_{\rho\rho}-{{1}\over{\rho}}q_{\rho}+{1\over{\rho^2}}q
-q\left(2|q|^2-4\int^{\infty}_{\rho}{{|q|^2(\tau,t)}\over{\tau}}d\tau
\right)-{ib\over{d-bt}}(q+\rho q_{\rho})=0. \label{qrho}
\end{eqnarray}
Then the following $({\wt p},{\wt q},{\wt r},{\wt u})$ giving by
the conformal transformation:
\begin{eqnarray}
{\wt p}(x,y,t)&=&
-{{ib{\bar z}}\over{2(a+bt)}}\nn\\
{\wt q}(x,y,t)&=&{{e^{-i{{b(x^2+y^2)}\over{4(a+bt)}}}}\over{a+bt}}
{q}(R,T)e^{-i\theta}\nn\\
{\wt r}(x,y,t)&=&-
{{e^{i{{b(x^2+y^2)}\over{4(a+bt)}}}}\over{(a+bt)}}{\bar q}(R,T)e^{-i\theta}\nn\\
{\wt
u}(x,y,t)&=&-{1\over{(a+bt)^2}}\left({|q(R,T)|^2}-2\int_R^{\infty}
{{|q(\tau,T)|^2}\over{\tau}}d\tau\right)+{{b^2z{\bar z}
}\over{2(a+bt)^2}} ,\nn \label{45qq}
\end{eqnarray}
where $R={{\rho}\over{a+bt}}$, $T={{c+dt}\over{a+bt}}$ and
$\left(\begin{array}{cc}
              a & b\\
              c& d\end{array} \right)\in SL_2({\bf R})$ (i.e.,
              $ad-bc=1$),
is a solution to the system (\ref{8}) with the restrictions
(\ref{q=r}).
\end{Lemma}

\nd {\bf Proof}. The proof is just a lengthy verification step by
step. For example, we have
\begin{eqnarray}
&&i{\wt q}_t-({\wt q}_{xx}+{\wt q}_{yy})+2{\wt u}{\wt
q}-2({\overline{\wt p}}{\wt q})_z+2{\wt p}{\wt q}_{\bar z}+4|{\wt
p}|^2{\wt q}\nn \\
 &&\quad={{e^{-i{{b(x^2+y^2)}\over{4(a+bt)}}}}\over{(a+bt)^3}}\left(
iq_T-q_{RR}-{1\over R}q_R+{1\over {R^2}}q-q(2|q|^2-
2\int^{\infty}_{R}{|q|^2\over{\tau}}d\tau)\right.\nn\\
&&\qquad\left.-{ib\over{d-bT}}q-
{ib\over{d-bT}}Rq_R\right)e^{-i\theta}\nn\\
&&\quad=0; \nn\\ &&i{\wt r}_t+({\wt r}_{xx}+{\wt r}_{yy})-2{\wt
u}{\wt r}-2({\overline{\wt p}}{\wt r})_z+2{\wt p}{\wt r}_{\bar
z}-4|{\wt
p}|^2{\wt r}\nn \\
 &&\quad={{e^{i{{b(x^2+y^2)}\over{4(a+bt)}}}}\over{(a+bt)^3}}\left(
-i{\bar q}_T-{\bar q}_{RR}-{1\over R}{\bar q}_R+{1\over
{R^2}}{\bar q}-{\bar q}(2|q|^2-
2\int^{\infty}_{R}{|q|^2\over{\tau}}d\tau)\right.\nn\\
&&\qquad\left.+{ib\over{d-bT}}{\bar q}+
{ib\over{d-bT}}R{\bar q}_R\right)e^{-i\theta}\nn\\
&&\quad=0; \nn
\end{eqnarray}
and so the other equations. $\Box$

\begin{Remark} Combining the gauge equivalent structure of the 1+2
dimensional Landau-Lifshitz equation (\ref{4}) displayed in
Theorem 1 and 2 with the conformal invariance displayed in Lemma
1, we may obtain a corresponding conformal invariant property to
the 1+2 dimensional isotropic Landau-Lifshitz equation (\ref{4}).
We call it the partial conformal invariant of Eq.(\ref{4}) since
it is only proved for some special case. However, it is very
interesting to see if the 1+2 dimensional Landau-Lifshitz equation
admits a similar conformal invariance in general.
\end{Remark}

Now we concentrate ourselves on characterizing some useful
analytic characters of the Cauchy problem of the following
equation:
\begin{eqnarray}\left\{
\begin{array}{c}
iq_t+q_{\rho\rho}+{1\over{\rho}}q_{\rho}-{1\over{\rho^2}}q
+q\bigg(2|q|^2-4\int^{\infty}_{\rho}{{|q|^2(\tau,t)}\over{\tau}}d\tau
\bigg)-{ib\over{d-bt}}\bigg(q+\rho q_{\rho}\bigg)=0\\
q(\rho,0)=q_0(\rho).
\qquad\qquad\qquad\qquad\qquad\qquad\qquad\qquad\qquad\qquad\qquad\quad
\end{array} \right.\label{qrho1}
\end{eqnarray}
or equivalently
\begin{eqnarray}\left\{
\begin{array}{c}
iQ_t+\Delta Q
+Q\bigg(2|Q|^2-4\int^{\infty}_{\rho}{{|Q|^2(\tau,t)}\over{\tau}}d\tau
\bigg)-{ib\over{d-bt}}\bigg(Q+\rho Q_{\rho}\bigg)=0\\
Q(\rho,0)=Q_0(\rho).
\qquad\qquad\qquad\qquad\qquad\qquad\qquad\qquad\qquad\quad
\end{array} \right.\label{qrho11}
\end{eqnarray}
where $Q(\rho,\theta,t)=q(\rho,t)e^{-i\theta}$ and
$Q_0(\rho,\theta)=q_0(\rho)e^{-i\theta}$. It is easy to see that
Eq.(\ref{qrho}) is complex conjugate to Eq.(\ref{qrho1}) and vice
versa. Moreover, if the initial data $q_0(\rho)$ satisfies
$q_0(\rho)e^{-i\theta}\in H^2({\bf R}^2)$, the local existence of
$H^2({\bf R}^2)$-solutions of the form $Q=q(\rho,t)e^{-i\theta}$
to the Cauchy problem (\ref{qrho11}) (or equivalently the local
existence of solutions $q(\rho,t)$ to the Cauchy problem
(\ref{qrho1}) such that $q(\rho)e^{-i\theta} \in H^2({\bf R}^2)$)
can be deduced directly from the known theory of nonlinear
Schr\"odinger-type equations or indirectly from the local
$H^3({\bf R}^2)$-existence and uniqueness result of the
Landau-Lifshitz equation proved by Sulem, Sulem and Bardos in
\cite{SuSu}, Theorem 1 and Lemma 1 (in this way, we require
additionally $\rho^2q_0(\rho)e^{-i\theta}\in L^2({\bf R}^2)$,
$\rho (q_0(\rho))_{\rho}e^{-i\theta}\in L^2({\bf R}^2)$ such that
$S_z\bigg|_{t=0}, S_{\bar z}\bigg|_{t=0}\in H^2({\bf R}^2)$ in the
case of $c=0$ in Lemma 1).

\begin{Lemma}
If $q(\rho,t)$ ($0\le t<T_0$ for some $0<T_0\le \infty$) is the
unique solution to the Cauchy problem (\ref{qrho1}) which
satisfies $Q(\rho,\theta,t)=q(\rho,t)e^{-i\theta}\in H^2({\bf
R}^2)$ (i.e., $Q(\rho,\theta,t)=q(\rho,t)e^{-i\theta}\in H^2({\bf
R}^2)$ solves the Cauchy problem (\ref{qrho11})), then

i).  For any $t$ with $0\le t<T_0$,\begin{eqnarray}
||q||_{L^2({\bf R}^2)}= ||q_0||_{L^2({\bf R}^2)} \quad
(\hbox{i.e.}\quad ||Q||_{L^2({\bf R}^2)}= ||Q_0||_{L^2({\bf
R}^2)}). \label{est1}
\end{eqnarray}

ii). There exists a positive constant $C$ such that, if $b<0$ and
$d>0$,
\begin{eqnarray}||\nabla^2Q||_{L^2({\bf R}^2)}^2\ge
{1\over{Ct+{1\over{||\nabla^2Q_0||_{L^2({\bf R}^2)}^2}}}}, \quad
0\le t<T_0. \label{est4}
\end{eqnarray}
\end{Lemma}
\nd{\bf Proof}. We multiply Eq.(\ref{qrho1}) by $2{\bar q}$ and
take the imaginary part of the result to get
$$
{{d}\over {dt}}|q|^2+2\nabla(\hbox{Im}{\bar q}\nabla
q)-{{b}\over{d-bt}}\bigg(2|q|^2+\rho|q|^2_{\rho}\bigg)=0.$$
Integrating the above equation on ${\bf R}^2$, we get
$${{d}\over
{dt}}\int_{{\bf R}^2}|q|^2dxdy=0.$$ Thus we have (i) by solving
this trivial ODE. In order to prove ii), we make the
transformation: $Q(\rho,\theta,t)\to {\wt
Q}(\rho,\theta,t)=e^{-Et}Q(\rho,\theta,t)$, where $E$ is a
positive constant which will be determined latter, to obtain
equivalently the following equation for ${\wt Q}$:
\begin{eqnarray}
i{\wt Q}_t+\Delta_{{\bf R}^2}{\wt Q} +{\wt Q}e^{2Et}\bigg(2|{\wt
Q}|^2-4\int^{\infty}_{\rho}{{|{\wt Q}|^2(\tau,t)}\over{\tau}}d\tau
\bigg)-{ib\over{d-bt}}\bigg({\wt Q}+\rho {\wt
Q}_{\rho}\bigg)+iE{\wt Q} =0\label{qrho2}
\end{eqnarray}
Taking the derivative with respective to $z$ to the both sides of
(\ref{qrho2}), we see
\begin{eqnarray}
&&i({\wt Q}_z)_t+\Delta_{{\bf R}^2}({\wt Q}_z) +4e^{2Et}{\wt Q}_z
|{\wt Q}|^2+2e^{2Et} {\wt Q}^2({\bar {\wt Q}})_z +4e^{2Et}{\wt
Q}_z\int^{\infty}_{\rho}{{|{\wt
Q}|^2(\tau,t)}\over{\tau}}d\tau\nn\\&&\qquad+ 4e^{2Et}{{{\wt
Q}|{\wt Q}|^2e^{-i\theta}}\over{\rho}} -{ib\over{d-bt}}\bigg(2{\wt
Q}_z+\rho ({\wt Q}_z)_{\rho}\bigg)+iE{\wt Q}_z=0.\label{qrho3}
\end{eqnarray}
Similarly we have
\begin{eqnarray}
&&i({\wt Q}_{\bar z})_t+\Delta_{{\bf R}^2}({\wt Q}_{\bar z})
+4e^{2Et}{\wt Q}_{\bar z} |{\wt Q}|^2+2e^{2Et} {\wt Q}^2({\bar
{\wt Q}})_{\bar z} +4e^{2Et}{\wt Q}_{\bar
z}\int^{\infty}_{\rho}{{|{\wt
Q}|^2(\tau,t)}\over{\tau}}d\tau\nn\\&&\qquad+ 4e^{2Et}{{{\wt
Q}|{\wt Q}|^2e^{i\theta}}\over{\rho}} -{ib\over{d-bt}}\bigg(2{\wt
Q}_{\bar z}+\rho ({\wt Q}_{\bar z})_{\rho}\bigg)+iE{\wt Q}_{\bar
z}=0.\label{qrho4}
\end{eqnarray}
We continue to take the derivative with respect to $z$ to the both
sides of (\ref{qrho3}) and to have
\begin{eqnarray}
&&i({\wt Q}_{zz})_t+\Delta_{{\bf R}^2}({\wt Q}_{zz}) +4e^{2Et}{\wt
Q}_{zz}
|{\wt Q}|^2+4e^{2Et}{\wt Q}_z^2{\bar {\wt Q}}+8e^{2Et}{\wt Q}|{\wt Q}|^2_z\nn\\
&&\qquad+2e^{2Et} {\wt Q}^2({\bar {\wt Q}})_{zz} -4e^{2Et}{\wt
Q}_{zz}\int^{\infty}_{\rho}{{|{\wt
Q}|^2(\tau,t)}\over{\tau}}d\tau+
8e^{2Et}{{{\wt Q}_z|{\wt Q}|^2e^{-i\theta}}\over{\rho}}\nn\\
&&\qquad+ 4e^{2Et}{\wt Q}\left(({{|{\wt
Q}|^2}\over{\rho}})_{\rho}- {{|{\wt
Q}|^2}\over{\rho^2}}\right)e^{-2i\theta}
-{ib\over{d-bt}}\bigg(3{\wt Q}_{zz}+\rho ({\wt
Q}_{zz})_{\rho}\bigg)+iE{\wt Q}_{zz}=0.\label{qrho5}
\end{eqnarray}
We multiply Eq.(\ref{qrho5}) by $2{\overline{{Q}_{zz}}}$, take the
imaginary part of the resulting expression and integrate it on
${\bf R}^2$ to have
\begin{eqnarray}
&&{{d}\over {dt}}\int_{{\bf R}^2}|{\wt Q}_{zz}|^2+E\int_{{\bf
R}^2}| {\wt Q}_{zz}|^2-4{b\over{d-bt}}\int_{{\bf R}^2}| {\wt
Q}_{zz}|^2=-8 e^{2Et}\int_{{\bf R}^2}\hbox{Im}\bigg({\wt
Q}_z^2{\bar
{\wt Q}}{\overline{{{\wt Q}}_{zz}}}\bigg)\nn\\
&&\qquad-16e^{2Et} \int_{{\bf R}^2}\hbox{Im}\bigg({\wt Q}_z{\wt
Q}{\bar {\wt Q}}_z{\overline{{{\wt Q}}_{zz}}}\bigg)-4e^{2Et}
\int_{{\bf R}^2}\hbox{Im}\bigg({{\wt Q}}^2({\bar
{\wt Q}})_{zz}{\overline{{{\wt Q}}_{zz}}}\bigg)\nn\\
&&\qquad-16e^{2Et} \int_{{\bf R}^2}\hbox{Im}\bigg({\wt Q}_z{{|{\wt
Q}|^2}\over{\rho}}e^{-i\theta}
{\overline{{{\wt Q}}_{zz}}}\bigg)\nn\\
&&\qquad-8e^{2Et} \int_{{\bf R}^2}\hbox{Im}\left({\wt Q}(|{\wt
Q}|^2_ze^{-i\theta}/\rho- 2|{\wt
Q}|^2e^{-2i\theta}/\rho^2)e^{-2\theta}{\overline{{{\wt
Q}}_{zz}}}\right).
 \label{c4}
\end{eqnarray}
Since we have the following estimates (here $C$ stands for
different constants):
\begin{eqnarray}
\left|\int_{{\bf R}^2}\hbox{Im}\bigg({{{\wt Q}|{\wt
Q}|^2e^{-2i\theta}}\over{\rho^2}}{\overline{{\wt
Q}_{zz}}}\bigg)\right|&\le& \left(\int_{{\bf R}^2}| {\wt
Q}_{zz}|^2\right)^{1/2}\left(\int_{{\bf
R}^2}{{| {\wt Q}|^6}\over{\rho^4}}\right)^{1/2}\nn\\
&\le&||\nabla^2 {\wt Q}||_{L^2}||{\wt
Q}||_{L^{\infty}}\left(\int_{{\bf R}^2}{{|
{\wt Q}|^4}\over{\rho^4}}\right)^{1/2}\nn\\
&\le&C||\nabla^2 {\wt Q}||_{L^2}||Q||_{L^{\infty}}\left(\int_{{\bf
R}^2}|
{\wt Q}_{\rho}|^4\right)^{1/2}\quad \hbox{(Hardy inequality)}\nn\\
&\le&C||\nabla^2 {\wt Q}||_{L^2}^2||{\wt
Q}||^2_{L^{\infty}},\label{e11}
\end{eqnarray}
where we have used the Gagliardo-Nirenberg inequality: $||\nabla
f||_{L^{4}}\le C ||f||_{L^{\infty}}^{1/2}||\nabla^2
f||_{L^2}^{1/2}$ for some constant $C$ in the last inequality,
\begin{eqnarray}
\left|\int_{{\bf R}^2}\hbox{Im}\bigg({\wt Q}_z{{|{\wt
Q}|^2}\over{\rho}}e^{-i\theta} {\overline{{\wt
Q}_{zz}}}\bigg)\right|&\le& \left(\int_{{\bf R}^2}| {\wt
Q}_{zz}|^2\right)^{1/2}\left(\int_{{\bf R}^2}{{| {\wt
Q}|^8}\over{\rho^4}}\right)^{1/4}\left(\int_{{\bf
R}^2}| {\wt Q}_z|^4|\right)^{1/4}\nn\\
&\le&||\nabla^2 {\wt Q}||_{L^2}||{\wt
Q}||_{L^{\infty}}\left(\int_{{\bf R}^2}{{| {\wt
Q}|^4}\over{\rho^4}}\right)^{1/4}\left(\int_{{\bf
R}^2}| {\wt Q}_z|^4|\right)^{1/4}\nn\\
%&\le&C||\nabla^2 Q||_{L^2}||Q||_{L^{\infty}}\left(\int_{{\bf
%R}^2}| Q_{\rho}|^4\right)^{1/4}\left(\int_{{\bf
%R}^2}| Q_z|^4|\right)^{1/4} \hbox{(Hardy)}\nn\\
&\le&C||\nabla^2 {\wt Q}||_{L^2}||{\wt
Q}||_{L^{\infty}}\left(\int_{{\bf R}^2}| {\wt
Q}_z|^4|\right)^{1/2}\quad\hbox{(Hardy)}
\nn\\
&\le&C||\nabla^2 {\wt Q}||_{L^2}^2||{\wt
Q}||^2_{L^{\infty}}\label{e12}
\end{eqnarray}
and
\begin{eqnarray}
||{\wt Q}||_{L^{\infty}}^2\le C(||{\wt Q}||_{L^2}^2+||\nabla^2{\wt
Q}||_{L^2}^2), \quad {\wt Q}\in H^2({\bf R}^2). \label{e13}
\end{eqnarray}
We substitute (\ref{e11},\ref{e12},\ref{e13}) and other easily
obtained estimates into (\ref{c4}) to have
\begin{eqnarray}
&&{{d}\over {dt}}\int_{{\bf R}^2}|{\wt
Q}_{zz}|^2+\left(E-4{b\over{d-bt}}\right)\int_{{\bf R}^2}| {\wt
Q}_{zz}|^2\nn\\&&\qquad\ge-C||{\wt Q}||^2_{L^2}||\nabla^2{\wt
Q}||^2_{L^2}\ge -Ce^{2Et}(e^{-2Et}||Q_0||^2_{L^2}+||\nabla^2{\wt
Q}||^2_{L^2})||\nabla^2{\wt Q}||_{L^2}^2
 \nn
\end{eqnarray}
for some constant $C$. Here we have used the fact: $||{\wt
Q}||_{L^2}=e^{-Et}||Q_0||_{L^2}$ from part i) of this lemma. In a
completely similar way, after either taking derivative with
respect to ${\bar z}$ to the both sides of (\ref{qrho3}) or with
respect to $z$ to both sides of (\ref{qrho4}), and did as above,
we still have
\begin{eqnarray}
&&{{d}\over {dt}}\int_{{\bf R}^2}|{\wt Q}_{z{\bar
z}}|^2+\left(E-4{b\over{d-bt}}\right)\int_{{\bf R}^2}| {\wt
Q}_{z{\bar z}}|^2\nn\\&&
\qquad\qquad\ge-Ce^{2Et}(e^{-2Et}||Q_0||^2_{L^2}+||\nabla^2{\wt
Q}||^2_{L^2})||\nabla^2{\wt Q}||_{L^2}^2 \nn \label{c6}
\end{eqnarray}
for some positive constant $C$. From the above two inequalities,
we obtain
\begin{eqnarray}
{{d}\over {dt}}\int_{{\bf R}^2}|\nabla^2 {\wt
Q}|^2+\left(E+C||Q_0||^2_{L^2}-4{b\over{d-bt}}\right)\int_{{\bf
R}^2}|\nabla^2 {\wt Q}|^2\ge -Ce^{2Et}||\nabla^2{\wt Q}||^4_{L^2}
 \label{c9}
\end{eqnarray}
for some positive constant $C$. Here we have used the fact that
the norm $\int_{{\bf R}^2}(|{\wt Q}_{zz}|^2+ |{\wt Q}_{z{\bar
z}}|^2)dxdy$ is equivalent to the norm $\int_{{\bf
R}^2}|\nabla^2{\wt Q}|^2dxdy$. That is to say, there is an
absolute positive constant $C_0$ such that $\int_{{\bf
R}^2}|\nabla^2{\wt Q}|^2dxdy\le C_0\int_{{\bf R}^2}(|{\wt
Q}_{zz}|^2+ |{\wt Q}_{z{\bar z}}|^2)dxdy\le {1\over
{C_0}}\int_{{\bf R}^2}|\nabla^2{\wt Q}|^2dxdy$. One may see, when
we set $E=C||Q_0||^2_{L^2}-4{b\over d}>0$, that the differential
inequality (\ref{c9}) implies
\begin{eqnarray}
{{d}\over {dt}}||\nabla^2 {\wt Q}||_{L^2}^2+2E ||\nabla^2{\wt
Q}||_{L^2}^2\ge -Ce^{2Et}||\nabla^2{\wt Q}||_{L^2}^4,
 \nn
\end{eqnarray}
which is equivalent to
\begin{eqnarray}
{{d}\over {dt}}{1\over{||\nabla^2 {\wt Q}||_{L^2}^2}}-2E
{1\over{||\nabla^2{\wt Q}||_{L^2}^2}}\le Ce^{2Et} \nn \label{c14}
\end{eqnarray}
or
\begin{eqnarray}
{{d}\over {dt}}\left({1\over{e^{2Et}||\nabla^2{\wt
Q}||_{L^2}^2}}\right)\le C.
 \nn\label{c15}
\end{eqnarray}
Therefore
\begin{eqnarray}
{1\over{e^{2Et}||\nabla^2 {\wt
Q}||_{L^2}^2}}-{1\over{||\nabla^2{\wt Q}_0||_{L^2}^2}}\le Ct,\quad
t\ge0.
 \nn
\end{eqnarray}
This shows (\ref{est4}) by substituting $||\nabla^2
Q||^2_{L^2}=e^{2Et}||\nabla^2{\wt Q}||^2_{L^2}$ and $||\nabla^2
Q_0||_{L^2}=||\nabla^2{\wt Q}_0||_{L^2}$. $\Box$

\bigskip

 We are in the position to prove our main result of this paper.
\begin{Theorem} There are $H^3({\bf R}^2)$-solutions
to the 1+2 dimensional  Landau-Lifshitz equation (\ref{4}), which
blow up in finite time.
\end{Theorem}
\nd{\bf Proof}. We only need to show the existence of some
solutions to the 1+2 dimensional isotropic Landau-Lifshitz
equation (\ref{4}) such that their $H^3({\bf R}^2)$-norms blow up
in finite time. For this purpose, we first take an initial
(complex) radial function $q_0(\rho)$ such that
$Q_0(\rho,\theta)=q_0(\rho)e^{-i\theta}\in H^2({\bf R}^2)$
%(and $\rho^2q_0(\rho)e^{-i\theta}\in L^2({\bf R}^2)$, $\rho
%(q_0(\rho))_{\rho}e^{-i\theta}\in L^2({\bf R}^2)$ in additional if
%necessary)
and $||Q_0||_{L^2({\bf R}^2)}$ is chosen to be so small that will
be specified below. Then, for any given $b, d$ with $b<0$ and
$d>0$, we solve the Cauchy problem (\ref{qrho1}) to get its unique
smooth solution $q(\rho,t)$ with $q(\rho,t)e^{-i\theta}\in
H^2({\bf R}^2)$. Setting $Q(\rho,\theta,t)= {\hat
q}(\rho,t)e^{-i\theta}$, where ${\hat q}(\rho,t)= {\bar q}(\rho)$,
we see that $Q(\rho,\theta,t)$ a smooth $H^2({\bf R}^2)$-solution
to (\ref{qrho}). There are only two possibilities to $Q$, say, a)
there is a finite $T_0>0$ such that $\lim_{t\to
T_0^-}||\nabla^2Q||_{L^2}=+\infty$ (short time $H^2$-existence)
and b) $||\nabla^2 Q||_{L^2}<\infty$ for any $0<t<\infty$ (long
time $H^2$-existence). In the followings we shall show separately
that, in either the case of the short or the long time existence,
there is a solution $S=S(z,{\bar z},t)$ to the Landau-Lifshitz
equation (\ref{5}) such that its $H^3({\bf R}^2)$-norm blows up in
finite time.

Before doing these, we first choose some real $a$ and $c$ such
that the $2\times2$ matrix $\left(\begin{array}{cc}
              a & b\\
              c& d\end{array} \right)\in SL(2,R)$, i.e. $ad-bc=1$.
The precise choice of $a$ and $c$ will be determined latter. From
Lemma 1, we see that
\begin{eqnarray}
{\wt p}(x,y,t)&=&
-{{ib{\bar z}}\over{2(a+bt)}}\nn\\
{\wt q}(x,y,t)&=&{{e^{-i{{b(x^2+y^2)}\over{4(a+bt)}}}}\over{a+bt}}
{\hat q}(R,T)e^{-i\theta}\nn\\
{\wt r}(x,y,t)&=&-
{{e^{i{{b(x^2+y^2)}\over{4(a+bt)}}}}\over{(a+bt)}}{\bar{\hat q}}(R,T)e^{-i\theta}\nn\\
{\wt u}(x,y,t)&=&-{1\over{(a+bt)^2}}\left({|{\hat
q}(R,T)|^2}-2\int_R^{\infty} {{|{\hat
q}(\tau,T)|^2}\over{\tau}}d\tau\right)+{{b^2z{\bar z}
}\over{2(a+bt)^2}} ,\nn
\end{eqnarray}
is a smooth solution to Eq.(\ref{8}) with the restriction
(\ref{q=r}), where  $R={{\rho}\over{a+bt}}$ and
$T={{c+dt}\over{a+bt}}$. By Theorem 1 and Theorem 2, there is a
smooth solution $S(z,{\bar z},t)$ to the Eq.(\ref{5}) which is
gauge equivalent to the solution $({\wt p},{\wt q},{\wt r},{\wt
u})$ of (\ref{8}) with the restrictions (\ref{q=r}). From the
formula: $S_{z}=2G U^{(\hbox{off-diag})}\sigma_3G^{-1}$, $S_{\bar
z}=-2G P^{(\hbox{off-diag})}\sigma_3G^{-1}$ deduced by the
relation (\ref{12}), where $U=\left(\begin{array}{cc}
  {\wt p} & {{\wt q}}\\
  {\wt r}& -{\wt p}
  \end{array}\right)$, $P=\left(\begin{array}{cc}
  \overline{{\wt p}} & \overline{{{\wt r}}}\\
  \overline{{\wt q}}& -\overline{{\wt p}}
  \end{array}\right)$ and the fact that $G\in
SU(2)$ is smooth (which implies that the absolute of every entry
of $G$ is not greater than 1), we have, as an entry, $|\wt q|\le
(|(s_1)_{z}|+|(s_2)_{z}|+|(s_3)_{z}|).$ Thus we obtain
\begin{eqnarray}
||S_{z}||_{L^2({\bf R}^2)}^2\ge  ||{Q}||^2_{L^2({\bf
R}^2)}=||Q_0||^2_{L^2({\bf R}^2)}.\nn\label{hao0}
\end{eqnarray}
In a similar way, we have the formulae:
\begin{eqnarray}
S_{zz}=2G\left(\begin{array}{cc}
              -2{\wt q}{\wt r} & 2{\wt p}{\wt q}-{{\wt q}_z}\\
              {2\wt p}{\wt r}+{{\wt r}_z}& 2{\wt q}{\wt r}\end{array}
              \right)G^{-1},\quad
S_{\bar z\bar z}=2G\left(\begin{array}{cc}
              -2\overline{{\wt q}}\overline{{\wt r}} &
              2\overline{{\wt p}}\overline{{\wt r}}+\overline{{\wt r}}_{\bar z}\\
              2\overline{{\wt p}}\overline{\wt q}-\overline{{\wt q}}_{\bar z}&
              2\overline{\wt q}\overline{\wt r}\end{array}
              \right)G^{-1}.\nn
\end{eqnarray}
%Thus $||S_{zz}||^2_{L^2({\bf R}^2)}\ge
%               ||2{\wt p}{\wt q}-{{\wt q}_z}||^2_{L^2({\bf
%R}^2)}={1\over{(a+bt)^2}}||{ Q}_Z||^2_{L^2({\bf R}^2)}$ and
%$||S_{\bar z\bar z}||^2_{L^2({\bf R}^2)}\ge
%               ||2{\wt p}{\wt q}-{\wt q}_{\bar z}||^2_{L^2({\bf
%R}^2)}={1\over{(a+bt)^2}}||{Q}_{\bar Z}||^2_{L^2({\bf R}^2)}$,
%where $Q=Q(R,\theta,T)=Q(Z,\bar Z,T)$. This shows that
%\begin{eqnarray}||S_{zz}||^2_{L^2({\bf
%R}^2)}+||S_{\bar z\bar z}||^2_{L^2({\bf R}^2)}\ge
%{1\over{(a+bt)^2}}||\nabla { Q}(T)||^2_{L^2({\bf R}^2)}.
%\label{hao1}
%\end{eqnarray}
\begin{eqnarray}
&&S_{zzz}=2G\left(\begin{array}{cc}
              -3({\wt q}{\wt r})_z & A\\
              B
              & 3({\wt q}{\wt r})_z\end{array}
              \right)G^{-1},\quad S_{\bar z\bar z\bar z}=2G\left(\begin{array}{cc}
              -3(\bar{\wt q}\bar{\wt r})_z & D\\
              E
              & 3(\bar{\wt q}\bar{\wt r})_z\end{array}
              \right)G^{-1},\nn\\
              &&
S_{zz\bar z}=2G\left(\begin{array}{cc}
              F
               & H\\
              J
              &-F \end{array}
              \right)G^{-1} \nn
\end{eqnarray}
where
\begin{eqnarray}
A&=&(2{\wt p}{\wt q}-{{\wt q}_z})_z -2{\wt p}(2{\wt p}{\wt
q}-{{\wt q}_z})-4{\wt q}^2{\wt r}\nn\\
&=&- {{e^{-i{{b(x^2+y^2)}\over{4(a+bt)}}}}\over{a+bt}}{
Q}_{zz}(R,\theta,T)-
4{{e^{-i{{b(x^2+y^2)}\over{4(a+bt)}}}}\over{(a+bt)^3}}|{Q}|^2{Q}\nn\\
B&=&({2\wt p}{\wt r}+{{\wt r}_z})_z+2{\wt p}({2\wt p}{\wt r}+{{\wt
r}_z})+4{\wt q}{\wt r}^2\nn\\
D&=&(2\bar{\wt p}\bar{\wt r}+{\bar{\wt r}_{\bar z}})_{\bar z}
+2\bar{\wt p}(2\bar{\wt p}\bar{\wt r}+{\bar{\wt r}_{\bar
z}})+4\bar{\wt r}^2\bar{\wt q}\nn\\
&=& {{e^{-i{{b(x^2+y^2)}\over{4(a+bt)}}}}\over{a+bt}}{Q}_{\bar
z\bar z}(R,\theta,T)+
4{{e^{-i{{b(x^2+y^2)}\over{4(a+bt)}}}}\over{(a+bt)^3}}|{Q}|^2\bar {Q} \nn\\
E&=&(2\bar{\wt p}\bar{\wt q}-{\bar{\wt q}_{\bar z}})_{\bar z}
-2\bar{\wt p}(\bar{2\wt p}\bar{\wt q}-{\bar{\wt q}_{\bar
z}})-4\bar{\wt q}^2\bar{\wt r} \nn\\
F&=&-2({\wt q}{\wt r})_{\bar z}+{\bar{\wt r}}({2\wt p}{\wt
r}+{{\wt r}_z})_z-\bar{\wt q}(2{\wt p}{\wt q}-{{\wt q}_z})\nn\\
H&=&(2{\wt p}{\wt q}-{{\wt q}_z})_{\bar z} +2\bar{\wt p}(2{\wt p}
{\wt q}-{{\wt q}_z})+4|\wt r|^2{\wt q}\nn\\
&=&-{{e^{-i{{b(x^2+y^2)}\over{4(a+bt)}}}}\over{a+bt}}{Q}_{z\bar
z}(R,\theta,T)+
4{{e^{-i{{b(x^2+y^2)}\over{4(a+bt)}}}}\over{(a+bt)^3}}|{Q}|^2{Q}\nn\\
J&=&({2\wt p}{\wt r}+{{\wt r}_z})_{\bar z}-2\bar{\wt p} ({2\wt
p}{\wt r}+{{\wt r}_z})-4|\wt q|^2{\wt r}.\nn
\end{eqnarray}
Thus we see from the above matrices equations that
\begin{eqnarray}
&&||S_{zzz}||_{L^2({\bf R}^2)}+||S_{\bar z\bar z\bar
z}||_{L^2({\bf R}^2)}+||S_{zz\bar z}||_{L^2({\bf
R}^2)}\nn\\
&&\qquad\ge ||A||_{L^2({\bf R}^2)}+||D||_{L^2({\bf R}^2)}
+||H||_{L^2({\bf R}^2)}\nn\\
%&&\qquad\ge C_0 \left({1\over{a+bt}}\right)^{4}\left(||\nabla^2 {
%Q}(T)+4|Q|^2{Q}(T)||^2_{L^2({\bf R}^2)}\right)
%\nn\\
&&\qquad\ge \left({1\over{a+bt}}\right)^{2}\left(C||\nabla^2 {
Q}(T)||_{L^2({\bf R}^2)}-48||{Q}(T)||^3_{L^6({\bf R}^2)}\right)
\label{hao2}
\end{eqnarray}
for some positive constant $C$. On the other hand, by using the
Sobolev inequality $||Q||_{L^6({\bf R}^2)}\le
C||Q||^{1/3}_{L^2({\bf R}^2)}||\nabla Q||^{2/3}_{L^2({\bf R}^2)}$
and the Gagliardo-Nirenberg inequality $||\nabla Q||_{L^2({\bf
R}^2)}$ $\le C||Q||^{1/2}_{L^2({\bf
R}^2)}||\nabla^2Q||^{1/2}_{L^2({\bf R}^2)}$, where $C$ stands for
different positive constants, we have $||Q||_{L^6({\bf R}^2)}^3$
$\le C_0||Q||^2_{L^2({\bf R}^2)}||\nabla^2 Q||_{L^2({\bf R}^2)}=
C_0||Q_0||^2_{L^2({\bf R}^2)}$ $||\nabla^2 Q||_{L^2({\bf R}^2)}$
for some constant $C_0$. Here we have used the conservation law i)
in Lemma 2. Substituting this inequality into (\ref{hao2}), we may
obtain
\begin{eqnarray}
&&||S_{zzz}(T)||_{L^2({\bf R}^2)}+||S_{\bar z\bar z\bar
z}(T)||_{L^2({\bf R}^2)}+||S_{zz\bar z}(T)||_{L^2({\bf R}^2)}\nn\\
&&\qquad\ge \left({1\over{a+bt}}\right)^{2}\left(C
-48C_0||Q_0||^2_{L^2({\bf R}^2)}\right)||\nabla^2Q(T)||_{L^2({\bf
R})}. \nn\label{hao01}
\end{eqnarray}
Therefore we may choose the initial data $Q_0$ in advance such
that $||Q_0||^2_{L^2({\bf R}^2)}$ is small enough that
$C-48C_0||Q_0||^2_{L^2({\bf R}^2)}>0$. Then we have
\begin{eqnarray}
||S_{zzz}(T)||_{L^2({\bf R}^2)}+||S_{\bar z\bar z\bar
z}(T)||_{L^2({\bf R}^2)}+||S_{zz\bar z}(T)||_{L^2({\bf R}^2)}\nn\\
\ge C \left({1\over{a+bt}}\right)^{2}||\nabla^2Q(T)||_{L^2({\bf
R}^2)} \label{hao02}
\end{eqnarray}
for some constant $C$ depending on  $||Q_0||_{L^2}$.

Now let's discuss the two different situations of the short or
long time existence mentioned before.

a) short time $H^2$-existence case. That is, there is a finite
time $T_0>0$ such that $\lim_{t\to
T_0^-}||\nabla^2Q||_{L^2}=+\infty$. Since $t={{aT-c}\over{d-bT}}$,
we may choose the entries $a,b,c,d$ with $ab<0$ in advance such
that $t_0={{aT_0-c}\over{d-bT_0}}>0$ and
$t_0={{aT_0-c}\over{d-bT_0}}<-{a\over b}$ (for example $a=1,b=-1$,
$d=1$ is small and $c=0$) since $-{a\over
b}-{{aT_0-c}\over{d-bT_0}}=-{{ad-bc}\over{b(d-bT_0)}}>0$. Thus we
see that there is an one to one corresponding between $t\in
[0,t_0)$ and $T\in [0,T_0)$ and
$$\lim_{t\to t_0^-}||\nabla^2{ Q}(T)||_{L^2}=\lim_{T\to
T_0^-}||\nabla^2{ Q}(T)||_{L^2}=+\infty.$$ Therefore, by noting
that $||S||_{H^3({\bf R}^2)}^2\ge C_0(||S_{zzz}||_{L^2({\bf
R}^2)}+||S_{zz\bar z}||_{L^2({\bf R}^2)}+||S_{\bar z\bar z\bar
z}||_{L^2({\bf R}^2)})^2+||S||_{H^2({\bf R}^2)}^2$ for some
positive constant $C_0$ and the estimate (\ref{hao02}), there is a
finite time $ t_0<-a/b$ such that
$$
\lim_{t\to t_0^-}||S||_{H^3({\bf R}^2)}\ge C_0 C
\left({1\over{a+bt_0}}\right)^{2}\lim_{T\to T_0^-}||\nabla^2{
Q}(T)||_{L^2}=+\infty.
$$
This shows that $||S||_{H^3({\bf R}^2)}$ blows up in this case of
short time $H^2$-existence of $Q$.

b) long time $H^2$-existence case. That is,
$||\nabla^2Q(T)||_{L^2}$ exists for any $T>0$. From (\ref{hao02}),
we see that
\begin{eqnarray} ||S||^2_{H^3({\bf
R}^2)}&\ge& C_0(||S_{zzz}(T)||_{L^2({\bf R}^2)}+||S_{\bar z\bar
z\bar
z}(T)||_{L^2({\bf R}^2)}+||S_{zz\bar z}(T)||_{L^2({\bf R}^2)})^2\nn\\
&\ge& C
\left({1\over{a+bt}}\right)^{4}||\nabla^2Q(T)||^2_{L^2({\bf R}^2)}
\label{hao4}
\end{eqnarray}
for some constant $C$ depending on $||Q_0||_{L^2({\bf R}^2)}$ when
$||Q_0||_{L^2({\bf R}^2)}$ is chosen suitably small. By lemma 2
ii), we have (if $b<0$ and $d>0$)
\begin{eqnarray}||\nabla^2 Q(T)||_{L^2({\bf R}^2)}^2\ge
{1\over{C_1T+{1\over{||\nabla^2Q_0||_{L^2({\bf R}^2)}^2}}}}, \quad
T>0. \nn
\end{eqnarray}
for some positive constant $C_1$ and hence
\begin{eqnarray}
||S||_{H^3({\bf R}^2)}^2 \ge C\left({1\over{a+bt}}\right)^{3}
{{(d-bT)}\over{C_1T+{1\over{||\nabla^2Q_0||_{L^2({\bf
R}^2)}^2}}}}. \label{hao5}
\end{eqnarray}
Here are have used the identity: $a+bt={1\over{d-bT}}$. Because
$T$ is an increasing function of $t$, when set the entries $a>0$,
$b<0$, $c=0$ and $d>0$ in advance, we see that there is one to one
corresponding between $t\in [0,-b/a)$ and $T\in [0,+\infty)$ and
$T\to +\infty$ as $t\to -b/a$. Thus, under this circumstance,
\begin{eqnarray}\lim_{t\to -b/a^-}\left({1\over{a+bt}}\right)^{3}
{{(d-bT)}\over{(C_1T+{1\over{||\nabla^2Q_0||_{L^2({\bf
R}^2)}^2}})}}=+\infty. \label{hao6}
\end{eqnarray}
Therefore, from (\ref{hao5}) and (\ref{hao6}), there is a finite
time $t_0\le -b/a$ such that
$$
\lim_{t\to t_0^-}||S||_{H^3({\bf R}^2)}=+\infty.
$$
This also shows that $||S||_{H^3({\bf R}^2)}$ blows up in finite
time in this case. The proof of Theorem 3 is completed. $\Box$

\bigskip

We would like to point out that, not like the nonlinear
Schr\"odinger equation with critical cases (\cite{Wei}), for
blowing-up $H^3({\bf R}^2)$-solutions to the 1+2 dimensional
Landau-Lifshitz equation (\ref{4}) constructed in the proof of
Theorem 3 b) one gets no their (up to the third derivatives)
point-wise blow-up information at the origin as $t\to -a/b$. So
our blowing up result does not contradict to the global existence
of $W^{m+1,6}({\bf R}^2)$-solutions to the Landau-Lifshitz
equations with small initial data due to Sulem, Sulem and Bardos
in \cite{SuSu}. Furthermore, we may see from Theorem 3 that it is
impossible to establish the global existence of $H^m({\bf
R}^2)$-solutions ($m\ge3$) to the Cauchy problem of the 1+2
dimensional Landau-Lifshitz equation for small initial data in
general. That is to say, the result of the global existence of
$W^{m,6}({\bf R}^2)$-solutions in \cite{SuSu} to the Cauchy
problem of the Landau-Lifshitz equation (\ref{4}) for small
initial data cannot be generalized to the case of energy estimates
in general. This indicates that the higher dimensional
Landau-Lifshitz equations may admit some unusual dynamical
properties.

In a similar way, the following conformal transformation of a
solution ${ Q}(\rho,t)$ to Eq.(\ref{QQ}):
\begin{eqnarray}
Q(\rho,t)\to {\wt {
Q}}(\rho,t)={{e^{-i{{b|x|^2}\over{4(a+bt)}}}}\over{a+bt}}{
Q}(R,T)\label{45}
\end{eqnarray}
is invariant, i.e., ${\wt {Q}}(\rho,t)$ is a solution to
Eq.(\ref{QQ}) too. The proof of this conclusion is a direct
computation and we omitted it here. It is well known that
nonlinear Schr\"odinger equations (\ref{eq}) with $\si=4/n+1$ have
localized finite energy solutions (\cite{Stra,Li,Wei}) which are
called solitary waves. Those are solutions of the form
$q(x;E_0)e^{iE_0t}$ with $E_0>0$, where $q(x;E_0)$ solves the
semi-linear elliptic equation $\Delta q-E_0q+|q|^{\si-1}q=0, q\in
H^1({\bf R}^n)$.  For our present Eq.(\ref{QQ}), the solutions of
the form
$$
\psi(\rho,t)=q(\rho,E)e^{iEt},
$$
where $E$ is a real constant and $q(\rho,E)$ solves
\begin{eqnarray}
q_{\rho\rho}+{1\over \rho}q_{\rho} -{1\over {\rho^2}}q-Eq +2|q|^2q
-4q\int_{\rho}^{\infty}{{|q|^2(\tau)}\over{\tau}}d\tau=0, \quad
q\in W^{1,\sigma}({\bf R}^2)\label{60}
\end{eqnarray}
for some $\sigma\ge2$, are also called solitary waves. We are
interested in $W^{1,\sigma}({\bf R}^2)$-solitary wave solutions
for general $p\ge2$ because of the $W^{1,\sigma}({\bf
R}^2)$-global existence of the Cauchy problem of the
Landau-Lifshitz equations (\ref{0}) with small initial data
obtained in \cite{SuSu} for $\sigma=6$ and in \cite{CSU} for
$\sigma=4$ and $n=2$. However, whether (\ref{60}) has a nontrivial
solution is unknown.

\bigskip

We finally remark that the system (\ref{8}) provides a new
mathematical point of view in investigating the 1+2 dimensional
isotropic Landau-Lifshitz equation (\ref{4}). The existence of
blowing-up $H^3({\bf R}^2)$-solutions to the Landau-Lifshitz
equation (\ref{4}) gives also an affirmative answer to the problem
proposed by Ding in \cite{Di} for Schr\"odinger maps.  However
whether there are $H^m({\bf R}^n)$-solutions to the
Landau-Lifshitz equation for general $n\ge3$ which blow up in
finite time is still unknown. We believe that some ideas displayed
in this paper will be helpful in understanding this problem.

\section * {\bf Acknowledgement}
The author is very grateful to Prof. Jiaxing Hong for his
suggestions and interested in this work. This work is partially
supported by the funds of doctoral bases of the Chinese Ministry
of Education, STCSM and NNFSC.

\end{document}